\documentclass{article}

\usepackage{natbib}

\usepackage{amssymb}

\usepackage[utf8]{inputenc}
\usepackage{graphicx,wrapfig}

\usepackage{endnotes}
\newenvironment{smquote}
{\begin{quote} \small}{\end{quote}}

\begin{document}
\title{Individual choice sequences 

History, development and use}  

\author{Joop Niekus\footnote{Institute for Language , Logic and Computation (ILLC),
University of Amsterdam,
e-mail j.m.niekus@uva.nl}}

\date{}

\maketitle

\begin{abstract}
  
We follow the history and development of Brouwer’s use of individual choice sequences up to the discovery of a method to apply them successfully in 1927. With the principles we derive from this first use we analyze in detail Brouwer’s work from that time onward. Our reconstruction uses only very basic principles. It  aligns exactly with Brouwer’s work after 1927 and, moreover, it gives a clear explanation of the proofs of his results and the terms he uses.
  
\end{abstract}

\bigskip

The central place in intuitionistic mathematics is occupied by the \emph{theory of continuum}. This theory is of a completely original character on account of the introduction of  an \emph{infinite sequence of arbitrary choices}.

Evert Beth in \cite{Beth1959}, p.~422.

\bigskip

\section{Introduction}
Brouwer introduced the real numbers by infinitely proceeding sequences with terms chosen more or less freely from mathematical objects previously acquired. This freedom of choice  can be restricted at any moment of the construction  of a sequence, up to complete determination. If the process of choice  is completely determined from the first term onward, we shall call the sequence \emph{lawlike}, otherwise we shall call it a \emph{choice sequence}\endnote{We followed Troelstra in restricting the expression choice sequences  for not completely determined sequences, e.g. \cite{Troelstra1977}. In the literature the expression choice sequence is sometimes used for any element of a spread, lawlike or not. But the distinction lawlike versus not completely determined has to be made in any case. Brouwer used in \cite{Brouwer1930A} in German `fertig' versus `unfertig', which made us use `complete' versus `incomplete,  in \cite{Niekus2010}. In his Post War articles Brouwer uses `sharp' for lawlike. We have chosen  the terminology from the standard texts.   }. Choice sequences are unique for Brouwer's intuitionism; understanding their nature seems crucial for understanding Brouwer's view of the continuum. 

But the concept is far from clear today. This is obvious from the reception of Brouwer's use of individual choice sequences, i.e.\ particular sequences not completely determined by a law. Their introduction in the late twenties was completely overlooked, and their post-war applications were misinterpreted and not recognized as choice sequences. Furthermore, the traditional theories for choice sequences are useless for handling the sequences Brouwer uses.

The traditional theories that have been developed to explain and formalize Brouwer's use of choice sequences, are based on Brouwer's work between 1918 and 1927. In that period Brouwer used only global properties of a continuum with choice sequences, such as the Continuity Principle CP, and not individual choice sequences, except for one dubious example at the end of the period, see below. The standard work on the subject is Anne Troelstra's \emph{Choice Sequences} (\cite{Troelstra1977}). A main result of this work is formed by the elimination theorems, according to which choice sequences are dispensable: what can be proved with choice sequences, can be proved without them. Still, in his final lectures Troelstra states that, on the basis of the elimination theorems, choice sequences have no mathematical value, but only philosophical.   For the standard theories substantial space has been reserved in the intuitionistic handbooks.

When Brouwer applied choice sequences after  the Second World War, he used a method that was commonly understood to be radically new. The sequences Brouwer defined were supposed to depend on the totality of the mathematical activity of an idealized mathematician, the IM. The method was controversial; it was claimed that Brouwer had lapsed into solipsism. Arend Heyting was even moderate with his comment on George Kreisel's presentation of a formalization of the method:

\begin{smquote}

It is true that Brouwer in his lecture   (`Conscious, philosophy and mathematics’, \cite{Brouwer1948C}) introduced an entirely new idea for which the subject is of essential importance. As I said, I feel that it is still questionable whether it is possible and whether it is good to introduce this idea in mathematics. [...]
In any case the method in question is not central in intuitionistic mathematics. [...] (\cite{Kreisel1967}, pp. 173).
\end{smquote}

Troelstra elaborated Kreisel's formalization in the standard-setting \cite{Troelstra1969}, resulting in the \emph{Theory of the Creative Subject}, the TCS. In this theory the particular sequences Brouwer defined were categorized as lawlike. From this assumption  Troelstra derived a paradox. In later expositions on Brouwer's \emph{method of the CS} the so called \emph{CS-sequences} are not always taken to be lawlike, but the TCS is still the standard theory for Brouwer's CS-arguments, see the recent \cite{Atten2018}. The TCS is controversial, also within intuitionism, and is treated as such in the last few pages of the intuitionistic handbooks.

In \cite{TroelstraDalen1988} Troelstra again discusses the paradox and concludes that the solutions he proposed in \cite{Troelstra1969} are not satisfactory, after which he just remarks that `In this connection the solution in \cite{Niekus1987} deserves further investigation'.  In our \cite{Niekus1987} we present a reconstruction without the dubious assumptions of the TCS. The defined sequences are sequences \emph{we} can construct, and for Brouwer's result we just need a principle for reasoning concerning the future. We claim that it is the application of a choice sequence that makes Brouwer's supposed new method special, and not the introduction of an IM.  

In our following papers \cite{Niekus2002} and \cite{Niekus2010}  we show that the method Brouwer applied was not new at all. We show that he uses it in his lectures from 1927 onward, as Brouwer remarked in \cite{Brouwer1948A}. The study of the texts of these lecturers fully supported our analysis of \cite{Niekus1987}. In all these papers the focus is on an instance of the example in \cite{Brouwer1948A}. 

In the present paper we follow a different route. We shall follow closely the history and development of Brouwer´s use of particular sequences, starting with his first unsuccessful attempt, followed by his breakthrough in the late twenties, up to his Post-War applications. It is especially the analysis of these Post-War examples that provides conclusive evidence for our position.

Brouwer's first use of a single choice sequence is in his proof of the negative continuity theorem in \cite{Brouwer1927B}. Brouwer uses a sequence with values chosen from a convergent sequence of reals and its limiting number, clearly not completely determined. The proof is not clear and allows different reconstructions. Moreover the result  is easily obtained when  a lawlike sequence is used, conforming to the elimination theorems. But for a long time this sequence has  been supposed to be the only particular sequence in the work of Brouwer. For that reason it has been reconstructed many times.  Brouwer's first attempt to apply a particular sequence may not have been a  successful one, but the applied technique in \cite{Brouwer1927B} extended with a new device, will  result in a method that turned out to be a breakthrough for Brouwer. 

The new device is the use of an assertion of which neither the absurdity nor the absurdity of its absurdity has been established. The idea of using such an assertion appeared for the first time in  \cite{Brouwer1926A}, where Brouwer gives an example of a virtual order $\prec$ which is not an order. Brouwer preferred this order relation in handling choice sequences over the natural order $<$.

Brouwer was well aware that he had made an important discovery. The new method would give direction to his further research of the continuum and he would use it in all of his future lectures: for the first time  in the Berlin Lectures of 1927, next in Vienna a year later and also in Geneva in 1934. The texts of these lectures support  the development a conception of choice sequence that arises naturally from Brouwer's constructive point of view. From this conception we shall derive the principles for handling Brouwer's  choice sequences.

We shall formulate the principles we derived from our conception of a choice sequence in the language of the TCS, but with a different interpretation of the basic term. In this interpretation two  axioms of the TCS are no longer valid, and neither is a principle derived from these axioms, Kripke's Schema KS, which is widely used for handling Brouwer's method. 

In the remainder of this paper we shall analyze several examples of Brouwer´s applications using the principles we adopt from our conception. Our first  example is from \cite{Brouwer1930A}, the text of the Vienna lecture; it is  the proof that  the intuitionistic continuum is not \emph{dense in itself}. In this proof a sequence is used with values chosen from a convergent sequence of reals   and its limiting number, just as in the \cite{Brouwer1927B} proof of the NCTH. These choices are made dependent on the results of the constructor of the sequence working on an unsolved problem, the assertion that determines the new method.

This combination of two techniques is refined further in the next examples we examine. They are from \cite{Brouwer1948C}, a Post-War paper in which Brouwer introduces the notion of a \emph{drift}, which enables him to reach his finest results. One may observe in our analysis how well the principles we have adopted are suited to explain Brouwer's arguments  and how well they can explain and express subtle distinctions Brouwer makes which cannot be made in the TCS.

The notion of drift appears again in  \cite{Brouwer1954A}, our final example. It is one of the rare occasions that Brouwer gives proofs of his results. One of these proofs   has been  used to support the standard theory TCS and Kripke's Schema KS. We claim that our analysis shows it does not. To the contrary, we claim the presentation of the proof is a conclusive argument against the TCS. We consider this analysis as an important result of the present paper.

\smallskip

The material of the paper is arranged as follows. We start in section 2 with the definition of a spread, the basic notion in intuitionism. We show Brouwer's proof of the existence of uncountable sets, using the continuity principle CP. After the introduction of the real numbers we show  the force of CP in proving that every real-valued function is continuous. After Brouwer's first restrictions on choice sequences in section 3, we discuss in section 4 Brouwer's proof of the NCTH in \cite{Brouwer1927A}, his first attempt to use a single choice sequence in a proof; the method was not yet successful. In section 5 we discuss the first appearance of the idea that would make the method work, the  use of an untested assertion. Brouwer was well aware that the new device made the method very powerful, see section 6. In section 7 we present the first appearance of the new method. From this example we derive the conception of a choice sequence, which enables us to derive principles for a formal treatment of the method. In section 8 we formulate these principles in the language of the TCS, resulting in the TICS. In section 9 our first reconstruction in the TICS with an example from the Vienna lecture of 1928 (\cite{Brouwer1930A}. In section 10 we express a distinction Brouwer makes in his Post-War papers in the TICS, and in section 11 we apply this distinction in the reconstruction of the examples of  \cite{Brouwer1948C}. In section 12 we reconstruct the examples of \cite{Brouwer1954A}, which turn out to be decisive arguments against the TCS. In section 13 we treat an attempt of Brouwer to eliminate the use of choice sequences, which is according to us unsuccessful. We end in the epilogue in section 14 with the conclusion that the TICS is superior to the TCS and the challenge for the traditional theories that is created by our results.

\section{The continuity principle CP} 

Brouwer started the intuitionistic reconstruction of mathematics in \cite{Brouwer1918B}. The article itself starts with the definition of a spread, the basic notion of intuitionism. Its importance for intuitionism can hardly be underestimated, as is demonstrated by Brouwer's introductory words for the spread definition in his Berlin Lectures of 1927:  

\begin{smquote}

   This set construction, which needs despite of its simple mental continence a lengthy description, and which carries alone the whole building of intuitionistic mathematics, is described as follows.\endnote{The original German text is: \begin{smquote} Diese Mengenkonstruktion, welche leider trotz ihres einfachen gedanklichen Inhaltes eine einigermassen langatmige Beschreibung erfordert, welche aber ganz allein das ganze Gebäude der intuitionistische Mathematik trägt, besteht in Folgendem: (\cite{Brouwer1991}, p.23).
\end{smquote} 

The translation to English here and elsewhere is ours, unless indicated differently.}

\end{smquote} 

Notwithstanding Brouwer's words, we shall give an adapted version of his definition:

\smallskip

A \emph{spread} is a \emph{law}, a \emph{building rule}, that regulates the construction of infinitely proceeding sequences, their terms chosen from a sequence $A$ of mathematical objects already constructed, the  \emph{founding sequence}. The building rule decides whether a choice from $A$ is admissible as the first term of a sequence under construction, and whether it is admissible as $n+1$-th term after a string of $n$ already admitted choices. After $n$ admitted choices there is always an admissible choice to continue the process. A thus constructed sequence, \emph{that in general has no complete description}, is called an element of the spread.

\smallskip

Before we focus on the emphasized phrase above, which is the subject of this paper, we first shall give a simple example  of a spread, which is the spread $C$. This spread has as its founding sequence the sequence of natural numbers $N$, and the building rule is that every choice is admissible, for the first term of a sequence as well as for every next. This spread is sometimes called the universal tree of natural numbers. Actually it is a description of the construction of that tree step by step. 

In his German texts Brouwer used for spread the German word for set \emph{Menge}. Closer to the classical notion of set is Brouwer's \emph{species} which is a property that can be possessed by a mathematical object. An object possessing the property is an element of the species. A spread is a species, but a species need not be a spread. Real numbers as introduced below are a species of species, not a spread.

An element of a spread is constructed term by term. If this process of construction has a complete description, i.e. if all of the values are  determined from the first term onward, the element is called  \emph{lawlike}. In case of the spread $C$ a lawlike element is given by a total function from $N$ to $N$. We shall call a sequence without such a complete description  a \emph{choice sequence}. But how is a choice sequence given to us? How do we use such incomplete objects? A lot of thought is spent on these questions in the secondary literature on intuitionism. But for Brouwer himself as well it was a life time struggle to clarify the notion of choice sequence. 

In the first years of his intuitionistic reconstruction Brouwer did not use specific properties of single choice sequences. The basic property he used is that in the most general case only an initial segment of an element of a spread is known. Brouwer uses this property to prove that the spread $C$ has a larger power than the natural numbers $N$ in the following way. 

\smallskip

\begin{smquote}

A function that assigns a natural number $n$ to an element of $C$ must fix that assignment on the basis  of an initial segment of that element. But then every element of $C$ sharing the same initial segment will be assigned the same natural number. So a $1-1$ function from $C$ to $N$ is impossible. Since a $1-1$ function from $N$ to $C$ is easily indicated, $C$ has a larger power than $N$ (\cite{Brouwer1918B}, p.13)\endnote{The original text, with $A$ denoting the natural numbers is:  `Die Menge $C$ ist gr\"{o}{\ss}er als die Menge $A$.  Ein Gesetz, das jedem Elemente $g$ von $C$  ein Element $h$ von $A$ zuordnet, muss n\"{a}mlich das Element $h$ vollst\"{a}ndig bestimmt haben nach dem Bekanntwerden eines gewissen Anfangssegmentes $\alpha$  der Folge von Ziffernkomplexen von $g$. Dann aber wird jedem Elemente von $C$, welches $\alpha$ als Anfangssegment besitzt, dasselbe Element $h$ von $A$ zugeordnet. Es ist mithin unm\"{o}glich, jedem Elemente von $C$ ein verschiedenes Element von $A$ zuzuordnen. Weil man anderseits in mannigfacher Weise jedem Elemente von $A$ ein verschiedenes Element von $C$ zuordnen kann, so ist hiermit der aufgestellte Satz bewiesen.' (\cite{Brouwer1918B}, p.~13.)}.
\end{smquote}

Brouwer applies here, without any explanation,   what has become  known as the \emph{continuity principle CP}. For    $\alpha$ and $\beta$ ranging over elements of $C$, for $f$ a function from $C$ to $N$ and for $\overline{\alpha_n}$ being (a code of) an initial segment of $\alpha$ of length $n$, we express CP in a formula by

\bigskip

CP ~~$\forall \alpha \exists n \forall \beta ( \overline{\alpha_n}=\overline{\beta_n} \rightarrow f(\alpha)=f(\beta))$

\bigskip

Classically CP is not valid: the function that assigns a $0$ to a sequence not containing a term with value $0$ and a 1  otherwise is well defined in classical mathematics, but intuitionistically it is not. CP is Brouwer's most characteristic tool and it makes intuitionistic analysis diverge essentially  from classical mathematics. After the introduction of the real numbers we shall show  the strength of CP in proving that every real valued function is continuous.

\bigskip

 We introduce the real numbers by  the spread RNG. The founding sequence is the sequence of integers; each integer is admitted as a first term, and $z$ is admitted as the successor $a_{n+1}$  of $a_1, a_2..., a_n$ if $z= 2a_n \vee z=2a_n+ 1 \vee z=2a_n+2$.
    
    To an element  $a=a_1, a_2, ...$ of RNG  corresponds  a shrinking sequence of  \emph{ $\lambda ^n$-intervals} of rational numbers  $\lambda_{a_1}^{1}, \lambda_{a_2}^{2},..., \lambda_{a_n}^{n}, ...$ defined by  $\lambda_{a_n}^{n}= [ \frac {a_n}{2^n}, \frac{a_{n }+2}{2^n}]$. Each  $\lambda ^n$ is contained in its predecessor by the definition of RNG,  and its length is $2^{-n+1}$.

  \smallskip
  
Brouwer called the elements of RNG \emph{points}. Two points $a$ and $b$ \emph{coincide}, $aRb$, if for all $n$ and $m$ the intersection of $\lambda_{a_n}^n$ and $\lambda_{b_m}^m$ is non-empty. This relation can be shown to be an equivalence relation. The equivalence classes,\emph{ point cores}, are Brouwer's real numbers, which we shall denote by $x, y, ... $.

\smallskip

 Given an element $a$ of a point core $x$ we can always construct for any $n$ an $a'$ such that $a' \in x$ and $a'$ is \emph{centered up to $n$}, i.e $\forall m<n~ (a'_m < a'_{m+1}/2 < a'_m+1 $), by defining $a'_k = a_k$ if $k \geq n$ and, by downward recursion from $n$,  $a'_k=(a'_{k+1}-1)/2 $ if $k<n$.

\smallskip

Brouwer had a preference to work with point cores; in this we follow him, conforming the quotations from his work we shall use below. But it is not necessary to proceed in this manner. Points by themselves are well suited for intuitionistic analysis, and will make the treatment  of real numbers less complicated, see the work of Wim Veldman (e.g. \cite{Veldman2000} or \cite{Veldman2021}). All handling of real numbers will be done via their generating points as in the following definitions of the natural order $<$ and of real valued functions. 

\bigskip

Let $a \in x$, $b \in y$, $r$ a rational number, and define 
\smallskip

$x<r$ if $\exists n ( \frac{a_n+2}{2^n}<r )$,

\smallskip

$x<y$ if $\exists n ( a_n+2 < b_n ) $,

\smallskip

$x \# y$ if $x<y \vee y<x$,

\smallskip

$\mid x-y \mid < r$ if $\exists n (\frac{\mid a_n- b_n \mid +2}{2^n} < r)$. \\

Consequently, if $a_m=b_m$ then   $\frac{\mid a_m- b_m \mid +2}{2^m} = 2^{-m+1} < 2^{-m+2} $,

 \smallskip

so $|x-y|<2^{-m+2}$.
                
\bigskip

A real-valued function $f$ is a law that assigns to a point core $x$ of its domain a point core $y$, notation $f(x)=y$; $f$ is  naturally given by a function $f'$ from RNG to RNG such that if $aRb$ than also $f'(a)Rf'(b)$.

\smallskip

For $a \in RNG$ and $n \in N$ let $\overline {f'(a)_n}$, as in the previous section,  be the initial segment of $f'(a)$ of length $n$. Since $\overline{f'(a)_n}$ is a finite sequence of integers, it can be coded uniquely by a natural number, so we may apply CP:

\smallskip

CPF~~~ $\forall a \forall m \exists n \forall b (\overline{a_n} = \overline{b_n} \rightarrow \overline{f'(a)_m} = \overline{f'(b)_m})$ ($\bigotimes$).

\smallskip

With CPF one can prove that every real valued function is continuous:

\smallskip

A function $f$ is \emph{continuous in x} if $\forall r \exists q_r \forall y \mid x-y \mid < q_r \rightarrow \mid f(x)-f(y) \mid < r $; 

$f$ is  \emph{continuous} if it is continuous in every point of its domain.

\bigskip

Theorem:  Every real-valued  function is continuous.

\bigskip

Proof: Let $f$ be a function generated by $f'$,  let $x$ be a pointcore, $a \in x$ 

and $m_0 \in N$.

\smallskip

By CPF  there is an $n_0$ such that 

\smallskip

$\forall b (\overline{a_{n_0}} = \overline{b_{n_0}} \rightarrow \overline{f'(a)_{m_0+2}} = \overline{f'(b)_{m_0+2}}$ ($\bigotimes$).

\smallskip

Without loss of generality we may assume that $a$ is \emph{centered up to $n_0+2$}.

\smallskip

Suppose $\mid x-y \mid < 2^{-n_0-2}$, $b \in y$. 

\smallskip

Since $a$ centered up to $n_0+2$, $\lambda_{a_{n_0}}^{n_0}$ contains $ \lambda _{b_{n_0+2}}^{n_0+2}$, so we can define

\smallskip

 $b'$ by $b'_m=a_m$  if $m  \leq n_0$, and continue $b'$ by   $b'_m=b_m$ for $m>n_0$.

\smallskip

Then $b' \in y$ and, since $\overline {a_{n_0}} = \overline{b'_{n_o}}$, $\bigotimes$ implies ~~ $\overline{f'(a)_{m_0+2}}=\overline{f'(b')_{m_0+2}}$.

\smallskip

So $\mid f(x)-f(y) \mid < 2^{-m_0}$,  $f$ is continuous in $x$.

\bigskip

Brouwer never proved this continuity theorem, the CTH. In \cite{Brouwer1927B} he proved a stronger theorem: every real valued function on [0,1] is uniformly continuous, the UCTH. For this purpose he first proved the fan theorem, which enabled him, instead of our starting point  $\forall a \forall m \exists n \forall b (\overline{a_n} = \overline{b_n} \rightarrow \overline{ f'(a)_m} = \overline{f'(b)_m} )$, to use $\forall m \exists n \forall a \forall b (\overline{a_n } = \overline{b_n} \rightarrow \overline{f'(a)_m} = \overline{f'(b)_m})$. From this point Brouwer's proof is straightforward,  similar to the proof above. But the proof of the fan theorem has caused much discussion, see section 5. 

But before Brouwer proved the UCTH he proved a theorem weaker than the continuity theorem, the negative continuity theorem, the NCTH. This proof also raised a lot of discussion, see  section 4.

\section{Restrictions on choice sequences}

During the first years of his intuitionistic reconstruction Brouwer uses global properties of a continuum with choice sequences, as the CP, without specifying their character further. But from 1925 on forward he irregularly added notes to the spread definition concerning the restrictions that can be imposed on a choice sequence during the construction. In \cite{Brouwer1925A} he appended to the phrase `a thus constructed sequence...':

\begin{smquote}
Including the feature of their freedom of continuation, which after each choice can be limited arbitrarily (possibly to being fully determined, but in any case according to a spread law).
\end{smquote}

With the text of Brouwer's Berlin Lectures of 1927 an extension of this note was found:

\begin{smquote}

The freedom to proceed with a choice sequence can after every choice arbitrarily be restricted (possibly in dependence on events in the world of  mathematical thought of the choosing person, imposed on the choosing person) (resulting e.g. in complete determination, or determination by a spread  law).\endnote{ The text is from \cite{Troelstra1982} p.~473. The author remarks that the note has been found in the Brouwer Archive together with the text of Brouwer's Berlin lecture from 1927, and that it has certainly been written before World War II.  Our quotation is Troelstra's translation of the orininal Dutch text, which is as follows: `De vrijheid van voortzetting van de betrokken keuzerij (eventueel in den kiezer opgelegde afhankelijkheid van gebeurtenissen in de wiskundige wereld van de kiezer) kan na iedere keuze willekeurig( b.v. tot volledige bepaaldheid, of ook volgens een spreidingswet) worden beperkt,...', 

The note is not included in \cite{Brouwer1991}, which contains the text of the Berlin lectures.}
\end{smquote}

The note is much longer. The part we have deleted concerns the possible restrictions on restrictions, restrictions on restrictions on restrictions, and so on. These higher order restrictions return in \cite{Brouwer1942A}. But in a note in \cite{Brouwer1981}, probably written around 1950, Brouwer expresses his doubts about their use, and finally, in \cite{Brouwer1952B}(p.142), he rejects the use of higher order restrictions completely. We just  mentioned the higher order restrictions and do not give the full text. They show Brouwer's struggle with the notion of choice sequence, but they are of no importance for our subject and neither have we found any trace of their influence on Brouwer's mathematical practice. 

That is different for the sentence in the citation that a choice may depend on the mathematical experience of the choosing person. It is directly connected with a discovery that would change the course of Brouwer's practice, which is the subject of this paper. The remark does not appear as a footnote in Brouwer's Post-War papers, but is standardly added to the definition of a choice sequence, see \cite{Brouwer1981}, \cite{Brouwer1952B} and \cite{Brouwer1954A}. From this last reference is the next quote (p.7). An arrow is a sequence  $p_1, p_2,..$ constructed according to the spread law.

\begin{smquote}

Furthermore, the freedom in the generation of an arrow, without being completely abolished,  may, at any $p_v$, undergo some restriction, and this restriction may be intensified at further $p_v$’s.
Finally, all these interventions, [...], may, at any stage, be made to depend on possible future mathematical experiences of the creating subject. 
\end{smquote}

The expression `choosing person' in the first quote of this section  is Troelstra's translation from the Dutch of the word 'kiezer'. In his Post-War papers Brouwer has switched to English, and as we may observe, now he denotes the choosing person, i.e. the constructor of the sequence, with the expression `creating subject', making explicit his philosophy that mathematics is a creation of the human individual. However, the expression `creating subject' was interpreted quite differently, as we shall see, leading to a misinterpretation. 

Because of the thoughts Brouwer gave to how a choice sequence is constructed, it may not be surprising that the next step is to incorporate choice sequences into his mathematical practice. In the next section we shall discuss Brouwer's first attempt.

\section{The negative continuity theorem NCTH}
In  \cite{Brouwer1927B} Brouwer proved one of his most discussed results, the uniform continuity theorem UCTH, see section 5 below. But first he proved, as an introduction, the negative continuity theorem, the NCTH. We do not treat the latter theorem because of the result, which follows immediately from the CTH proved above, but because of Brouwer´s proof. In this proof Brouwer uses for the first time  a sequence which is clearly not completely determined, a choice sequence. 

It is his first attempt, and we shall see the applied method needed an extra device to make it work. But for a long time this sequence was supposed to be the only individual choice sequence in the work of Brouwer. It is the only one discussed in \cite{Troelstra1982}, a study in `the origin and development of choice sequences' and in \cite{Atten2007}, a study `on the Phenomenology of Choice Sequences'. 

For this reason, and because of the fact that the proof is not clear, the proof has been reconstructed many times. We shall discuss some of the reconstructions after Brouwer's proof below.

\bigskip

    A function is \emph{negatively continuous} in $x$ if  for some infinitely proceeding  sequence $x_1,~x_2,~...$ converging to $x$ 

  \smallskip

$ \forall r>0~ \neg \forall n \exists m \mid f(x_{m+n}) - f(x) \mid > r$.

\smallskip

A function $f$ is negatively continuous on its domain if it is negatively continuous in all of its points.

\bigskip

Theorem: Every full function is negatively continuous.

\bigskip

 Brouwer´s proof, adapted version:\endnote{ Brouwer's proof:
\begin{smquote}
Sei $f(x)$ eine volle Funktion, $\xi_0$ eine beliebieger Punktkern $x$ und $\xi_1, \xi_2,\ldots$ eine gegen $\xi_0$ positiv konvergierende Fundamentalreihe von Punktkernen $x$. Wir nemen nun einen Augenblick an, dass es eine natürliche Zahl $p$ und eine Fundamentalreihe von beständig wachsenden natürlichen Zahlen $p_1, p_2,\ldots$ gäbe, so dass $\mid f(\xi_{p_v}-f(\xi_0 \mid > \frac{1}{p}$ für jedes $v$, und definieren einen Punktkern $\xi_{\omega}$ des Einheitskontinuums, indem wir von einer unbegrenzten Folge $F_1$  von erzeugenden Intervallen eines zu $\xi_0$ gehörenden Punktes ausgehen, und sodann durch eine unbegrenzte Folge von Wahlen von $\lambda$-Intervallen in solcher Weise einen Punkt $F_2$ des Einheitskontinuums konstruieren, dass wir einstweilen für jede schon ins Auge gefasste natürliche Zahl $n$ die ersten $n$ Intervalle mit den ersten $n$ Intervallen von $F_1$ identisch wählen, uns aber die Freiheit vorbehalten, jederzeit, nachdem das erste, zweite, …, $(m-1)$-te und $m$-te Intervall gewählt worden sind, die Wahl aller weiteren (d.h.des $(m+1)$-ten, $(m+2)$-ten uns. ) Intervalle in der Weise festzulegen, dass entweder ein zu $\xi_0$ oder rein zu einem gewissen $\xi_{p_v}$ gehöriger Punkt erzeugt wird. Alsdann ist für den $F_2$ enthaltenden Punktkern $\xi_{\omega}$ die Funktion $f(x)$ nicht definiert, womit wir zu einer widerspruch gelangt sind, und unsere Annahme sich als unstatthaft erwiesen hat. Dies aber besagt die negative Stetigkeit der Funktion $f(x)$ (\cite{Brouwer1927B}, p.62).
\end{smquote}
}

Let $f(x)$ be a full function, $\xi_0$ a point core and $\xi _1, \xi _2, \ldots$ a sequence of point cores converging to $\xi_0$. Suppose a natural number $p$ and a sequence of steadily growing natural numbers $p_1, p_2,\ldots$ exist, such that $\mid f(\xi_{p_v})-f(\xi_0 )\mid > \frac{1}{p}$ holds for every $v$. Define a point core $\xi_{\omega}$ starting with an infinitely preceding sequence $F_1$ of $\lambda$-intervals generating $\xi_0$, and construct an infinitely proceeding sequence of $\lambda$-intervals  $F_2$, defining a point, as follows: for every considered natural number $n$ we let the first $n$ intervals be identical with the first  $n$ intervals of $F_1$, but we maintain the freedom, after the choice of the first, second, ... , $(m-1)$-th and $m$-th interval have been made, to determine the choice of all further intervals (i.e. the $(m+1)$-th, $(m+2)$-th etc. ) in such a way, that either the resulting point generates $\xi_0$ or it generates a certain  $\xi_{p_v}$. But then the function $f(x)$ is not defined for the point core $\xi_{\omega}$ defined by $F_2$, which is a contradiction. This proves the negative continuity of the function   $f(x)$.

\bigskip

 After the proof Brouwer remarked that `By this theorem, which is an immediate consequence of the intuitionistic point of view, and which I have applied often in talks and lectures since 1918, the supposition of the validity of the much stronger theorem 3 [the UCTH, j.n.] could become close'.  \endnote{ The original German text of this remark is `Durch dieses eine unmittelbare Konsequenz des intuitionistischen Standpunktes bildende, seit 1918  von mir vielfach in Vorlesungen und Gesprächen erwähnte Theorem 1 wird die Vermutung der Gültigkeit des unterstehenden viel mehr aussagenden Theorems 3 nahegelegt..' (\cite{Brouwer1927B}, p.62). }
 
 \smallskip

His claim that he used the proof from 1918 on in his lectures is enigmatic. Brouwer mentioned negative continuity once in \cite{Brouwer1923B}, but there is no trace  of this theorem before \cite{Brouwer1927B} in the work of Brouwer. Van Dalen suggests that mentioning 1918 could have to do with Brouwer claiming to have continuity already in 1918. This might explain the formulation of the NCTH. Fullness is not needed in the proof but formulated in this way the NCTH is classically equivalent with the UCTH.  

\smallskip

  It is the definition of the sequence $\xi_{\omega}$ that makes the proof special. The values of $F_2$, the sequence generating  $\xi_{\omega}$, are chosen from $F_1$ for every considered $n$, but one can never exclude the possibility of continuing as a sequence generating some  $\xi_{p_v}$. Which means that we never reach a value for  $F_2$ after which all other values are fixed. Brouwer uses for the first time  a particular choice sequence.

 In intuitionistic research there is no agreement about what Brouwer means with `an immediate consequence of the intuitionistic point of view'. Various researchers have given different interpretations of the proof. We shall give two characteristic ones.

 A slight modification of our proof of the CTH in the previous section gives the result immediately. Such kind of proof is Troelstra's first and favorite reconstruction in \cite{Troelstra1982}. He constructs a sequence by taking the values of a $\xi _0$ generating sequence, until the function value is  narrowed enough to draw a contradiction. So   Troelstra assumes that we can narrow the function value by taking enough steps in the generating sequence, which implies a local form of CP. A drawback of the reconstruction is that he could have proved continuity of the function with the same effort; it would be a proof similar to the one in the previous section. 

One might wonder whether CP is necessary for the proof. In  \cite{Veldman1982} Wim Veldman gives a reconstruction without the use of CP. It runs as follows: 

\smallskip

Let $k_7$ be such that it is the number of the first digit in the decimal expansion of $\pi$ were a row of ten $7$'s starts, notations as in Brouwer´s proof. Define a sequence  $F_2$ of $\lambda$-intervals as follows: as long as $n<k_7$ we choose the $n$-the term of $F_2$ to be equal to the $n$-the term of $F_1$.  But for $n\geq k_7$ we choose the terms of $F_2$ such, that $F_2$ generates $p_{k_7}$. Let $x_\omega$ be the real number generated by $F_2$. Since the existence of $k_7$  is not known, neither the impossibility of its existence, we cannot calculate $f(x_\omega)$, it is not defined.
\smallskip

Veldman does not conclude to a real contradiction, as Brouwer does. He interprets the `not' in the weak sense, that is $ \forall n \exists m \mid f(x_{m+n}) - f(x) \mid > r$ is just not valid. Is CP necessary for a real contradiction? Furthermore, the sequence Veldman uses is characteristic for Brouwer's intuitionism, we will meet such sequences later in this paper. But it is lawlike, and the sequence used by Brouwer clearly is not.

The problem of deriving a strong negation from the weak `not' without CP is solved in Posy's (\cite{Posy1976}) and Martino's (\cite{Martino1985})  by the application of the Theory of the Creative Subject, the TCS. This theory has been  developed for the reconstruction of a Brouwerian method after World-War II. It certainly makes sense to connect Brouwer's CS-arguments with the application of a choice sequence, but Brouwer's concept of a choice sequence had not yet crystallized out, and as we shall show in section 8, there is no base for the TCS in intuitionistic reasoning.

Our conclusion is that for the reconstruction of Brouwer's proof there seem to be two options. One assumes (some form of) CP and one can prove continuity with the same effort, see also the recent \cite{Fletcher2020}, or one reaches a weak negation as in Veldman's reconstruction above.

We think that Brouwer is not careful  in his statements  and that he intended to prove a weak negation. The use of CP at this place in \cite{Brouwer1927B} is not obvious. Brouwer used the CP in \cite{Brouwer1918B} without further explanation, but in \cite{Brouwer1927B} he formulated it after the proof of the NCTH. And with CP the proof is much like the last part of the proof of the UCTH, and too sophisticated to be called an immediate consequence of the intuitionistic point of view. 

Further, a choice sequence is not needed  at all. The result, weak or strong, can easily be reached  with a lawlike sequence. Actually, Brouwer's definition of a choice sequence as a description of a construction has a paradoxical flavor, in always choosing values from $\xi _0$, but keeping open the possibility of switching to a $\xi_n$. We can only conclude that Brouwer's attempt to support the UCTH by a proof of the NTCH from really basic intuitionistic principles has to be considered a failure. Unjustifiably one tried to find the thoughts behind choice sequences in the many reconstructions. 

Our claim that Brouwer intended to give a weak result, and that his use of a choice sequence is not successful and not needed at all, is supported by the other appearances of the NCTH we found in Brouwer's work. Twice again he proved the theorem, both labeled as an immediate consequence of the intuitionistic point of view, in contrast with the UCTH, which needs the fan theorem for its proof.  

In his 1927 Berlin Lectures he gives again the NCTH as a preparation for the UCTH. The sequence he uses is lawlike, and the result is that if the NCTH holds, a principle classically valid and intuitionisticly not, would also be valid, which means that Brouwer proves a weak result (\cite{Brouwer1991} p.50). But at that point  Brouwer had found a more successful way to apply single choice sequences. In the same Berlin Lectures he  presents this method for the first time, and it would give direction to all of his further work. 

The second appearance is in the last publication of Brouwer's work, which is the text of Brouwer´s Cambridge Lectures. In these lectures he uses the new method in an extensive way. But in his proof of the NCTH he uses a lawlike sequence, like Veldman as above (\cite{Brouwer1981}, p.81).

The new method is a combination of choosing values from a convergent sequence of reals, as in the proof of the NCTH discussed above, extended with a new device, which is making the choices dependent on the results of the constructor of the sequence working on a mathematical problem. For the results obtained by this method the use of choice sequences seems essential. Somewhat surprisingly, the new device is of no help in Brouwer's proof of the NCTH; only a weak result can be obtained when it is applied.

The idea of that device appeared, as far as we know, for the first time in 1926 when Brouwer introduced the virtual order, see next section.

\section{The virtual order}

 In  \cite{Brouwer1926A} Brouwer defines the two-place relation \emph{virtual order} $\prec$, see the conditions below. As mentioned above in this paper it is the first time that the idea appears that Brouwer uses for his new device.

 The relation $x\prec y$ defined by $(x \prec y) \Leftrightarrow (x\neq y \wedge \neg y<x)$, with $<$ the natural order and $x$ and $y$ real numbers, fulfills the conditions of the virtual order. This order $\prec$ proved to be very appropriate while handling single choice sequences.   From the moment Brouwer introduced them, which was in the Berlin Lectures of 1927, he used this virtual order in his analysis of the continuum, and he would do so the rest of his career.

 None of the researchers after Brouwer followed him in using the negative property $\prec$ in stead of the natural $<$. We will follow the practice of these researchers in the next sections, except when discussing Brouwer's writings closely.
 
 \bigskip

A species is virtually ordered by the two place relation $\prec$ if it satisfies the following conditions. 

\smallskip

1. The relations $r=s$, $r \prec s$ and $r \succ s$ exclude each other.

\smallskip

2. From $r=u$, $s=v$ and $r \prec s$ it follows that $u \prec v$.

\smallskip
3. If both $r \succ s$ and $r=s$ are impossible, it follows that $r \prec s$. 

\smallskip
4. If both $r \succ s$ and $r \prec s$ are impossible, it follows that $r=s$. 

\smallskip
5. From $r \prec s$ and $s \prec t$ it follows that $r \prec t$. 

\bigskip

 In \cite{Brouwer1927C}   Brouwer proves that a relation is a virtual ordering if and only if    it is  an \emph{inextensible ordering}, i.e. if $a \prec b$ would be consistent with the ordering, then $a \prec b$ is already the case.

After the definition of virtual order Brouwer gives the following distinctions, which will appear further in this paper:

\smallskip
A virtually ordered species is \emph{ordered} if there is an ordering relation for any two different elements of $P$, i.e.

$\forall xy(x\neq y \rightarrow (x\prec y \vee y \prec x))$.

\smallskip
A virtual order is a \emph{complete order} if it is ordered and discrete, i.e.

$\forall xy(x=y \vee x \prec y \vee y \prec x)$.

\smallskip

In a footnote Brouwer enlightens these definitions with an example. In this example Brouwer introduced the idea that would give him the device for successfully handling choice sequences.

\begin{smquote}
That even a virtual ordered species $P$ consisting of two different elements $a$ and $b$ does not have to be ordered is shown by the following example: Let $d_v$ be the $v$-the cipher in the decimal expansion of $\pi$ and let $m=k_1$ if in the decimal expansion of $\pi$ at  $d_m d_{m+1} \ldots d_{m+9}$    for the first time in the expansion of $\pi$ the sequence $012345679$ occurs. Further, let $a \succ b$  when the existence or the absurdity of the absurdity of the existence of $k_1$ holds, and $a \prec b$ when the absurdity of the existence of $k_1$ holds. These definitions fulfill the 5 ordering principles, so they define a virtual ordering on $P$. But it is not decided that the tuple $(a,b)$ is an ordered pair, so the species $P$ is, with this definition,  not an ordered one (note  in \cite{Brouwer1926A}, p.~455).\endnote{Da{\ss} schon ein virtuell geordnete Spezies $P$ von zwei verschiedenen Elementen $a$ und $b$ nicht notwendig geordnet zu sein braucht, geht aus folgendem Beispiel hervor: Es sei $d_v$ die $v$-te Ziffer hinter dem Komma der Dezimalbruchentwickelung von $\pi$    und  $m=k_1$, wenn es sich in der vortschreitenden Dezimalbruch- entwickelung von $\pi$ bei $d_m$ zum ersten male ereignet, da{\ss} der Teil $d_m d_{m+1}...d_{m+9}$ dieser Dezimalbruchentwickelung eine Sequenz 0123456789 bildet. Es sei weiter $a>b$, wenn sich von $k_1$ entweder die existenz oder die Absurdit\"{a}t der Absurdit\"{a}t der Existenz, und $a<b$, wenn sich von $k_1$ die Absurdit\"{a}t der Existez herleiten l\"{a}{\ss}t. Diese Festsetzung gen\"{u}gt den f\"{u}nf Ordnungsbeziehungen, bestimmt also eine viruelle Ordnung von $P$. Weil aber nicht entschieden ist, da{\ss} da{\ss} f\"{u}r das Elementepaar $(a,b)$ eine ordnende Relation besteht, so liegt keine Ordnung der Spezies $P$ for.  
\cite{Brouwer1975}, p.323.}
\end{smquote}

We can work out this  example as follows:

Let $\varphi$ be the assertion $\exists n( k_1=n)$. Then Brouwer defines 

\smallskip

(1) $(\varphi \vee \neg \neg \varphi) \leftrightarrow a \succ b$, 

\smallskip

(it is curious that Brouwer does not simplify $\varphi \vee \neg \neg \varphi$ to $\neg \neg \varphi$), and

\smallskip

(2) $\neg \varphi \leftrightarrow a \prec b$ .

\smallskip

To show that a virtual order has been  defined we have to check the conditions 

1. up to 5. The only non-trivial case to check is 3.:

\smallskip

Suppose $\neg b \succ a$. 

Because of (2) we would have $\neg \neg \varphi$, so because of (1): $a \succ b$.

We cannot conclude to $\varphi$ from $\neg b \succ a$, so $\varphi \leftrightarrow a \succ b$ would not be sufficient 

to make it a virtual order.

And finally, as long as neither $\neg \varphi$ nor $\neg \neg \varphi$ has been established, the species 

$P$  is not ordered by $\prec$.

\smallskip

The property that neither the absurdity of the existence of $k_1$, nor the absurdity of the absurdity of its existence has been established, is an example of an assertion that has not been \emph{tested}. Remark that if one of these possibilities is established, the species above becomes completely ordered by the thus defined relation $\prec$. We will see that connecting the choices from a convergent sequence of reals with the results working on an untested assertion provides Brouwer the method for a successful application of single choice sequences.

\section{Nothing new?}

After the negative continuity theorem Brouwer proved his most famous and most discussed result: the uniform continuity theorem, the UCTH. The  starting point of the proof is that a real valued function has the properties we used in our proof of the CTH section 2, i.e. it  is given by a function from RNG to RNG, on which the continuity principle CP is applicable. Given such a function, Brouwer reasons, we also have a proof $h$ that the function has those properties. Analyzing $h$ Brouwer deduces a downward recursion principle for trees, named Bar Induction. With Bar Induction he proves straightforwardly the Fan Theorem, called Fundamental Theorem by Brouwer, and finally, using the Fan Theorem, uniform continuity, see section 2. The assumptions about $h$ leading to Bar Induction are subject of research up to the present day. Extensive analysis are e.g in   \cite{KleeneVesley1965}, \cite{Dummett1977} and \cite{Atten2004}. A recent contribution is \cite{Veldman2021}.

It is the proof of Bar Induction that raises questions, not the principle itself. It is generally accepted as a reasonable principle; it appears in intuitionistic formal   systems as an axiom. Brouwer was proud of the result, but not satisfied with the proof, which is demonstrated by a note found with the text of a 1952 lecture, presented by Walter van Stigt in his \emph{Brouwer´s Intuitionism} (\cite{Stigt1990} : 

\begin{smquote}

For Fans a truly wonderful theorem holds whose importance would justify to call it the  Fundamental Theorem of Intuitionism, but whose absolutely rigorous proof till now has not been sufficiently simplified(\cite{Stigt1990} p.93).
\end{smquote}

It made van Stigt wonder whether this dissatisfaction was the cause for the `abrupt ending of his intuitionistic reconstruction program', a judgement he based on a thorough study of Brouwer´s work after \cite{Brouwer1927B}. He stated:

\begin{smquote}

It remains an open question whether dissatisfaction with his proof of the Bar Theorem, his failure to find a `simple’, constructive alternative is the main cause of the demise of Brouwer’s program of re-constructing mathematics. It is a fact, however, that all his attempts at a systematic construction of an Intuitionistic Analysis abruptly end at this point.
The abrupt ending of Brouwer’s re-construction program at this stage is even more evident in his unfinished, unpublished books.  (\cite{Stigt1990}, p.93)
\end{smquote}

It is true that Brouwer's systematic reconstruction ends at this time. But it is not true that Brouwer did nothing new after \cite{Brouwer1927B}. One of his most remarkable achievements still had to come. This achievement was not a reconstruction result. It had to do with the basic subject of intuitionism, which is the continuum.  And apparently Brouwer was fully aware that he had found something special, see the next quote:

\begin{smquote}

There is a letter from the significist Henri Borel to Brouwer, right after the Berlin Lectures, in which he mentioned that "Gutkind wrote me that `you had made an enormous discovery', which `had to do with the foundations of logic'"  (\cite{Dalen1999A}, pp.394-395).
\end{smquote}

In \cite{Dalen2005} (p.553) van Dalen precedes the same quote with `There is a puzzling letter from Henri Borel...'. And after the quote he states `Neither in Brouwer's papers, nor his private notes is such a discovery recorded, assuming that Gutkind referred to a recent event.' After considering some possibilities, neither of them judged to be satisfactory, he concludes that `It is just possible that Brouwer had discussed the idea and the use of the creating subject as a tool to study the full continuum'. The author seems to be a bit reluctant to connect the `enormous discovery' with a method that has in intuitionistic research the status of obscure and controversial. But this status is, as we shall show,  based on a misinterpretation of the method. Brouwer had made a breakthrough in handling choice sequences. 

\smallskip
The letter mentioned in the quote is dated 1927, March 26. Gutkind was a friend of Brouwer's. Brouwer stayed in his house during the Berlin lectures of January 1927. In these lectures Brouwer introduced a method in which he applied particular choice sequences. It must be the case that Brouwer was referring to this method  in speaking of an enormous discovery.

\section{Berlin 1927 - What is a choice sequence?}
In 1927 Brouwer lectured in Berlin. The structure of these lectures would be a model for all his following lectures. He starts with a historical overview and  next he introduces the intuitionistic fundamentals  as spread,  species and the continuum. Thereafter he treats the notion of order,  and gives a further analysis of the continuum, and finally his proof of the UCTH. 

With his fundamental concepts Brouwer introduces a new distinction. Real numbers generated by lawlike sequences form the \emph{reduced continuum}, all sequences  of the real number generating spread, lawlike or not, generate the \emph{(full) continuum}. Moreover  Brouwer had found a method to explore this distinction. He had discovered a way to define sequences not completely determined by a law, choice sequences, and a method to apply these sequences in his study of the notion of order and in his analysis of the continuum. The method relied on combining the use of an untested proposition, as in the previous section, with the mathematician at work. 

Below we shall give the first occurrence of the new method, which is Brouwer´s proof that the full continuum is \emph{not ordered} by the natural $<$. First we shall show Brouwer's proof that the natural order $<$ is \emph{not a complete order} on the reduced continuum, i.e.\ $~x>0~ \vee ~ x=0 ~\vee x<0$ does not hold generally for lawlike real numbers:

\begin{smquote}
Further, we denote with $K_1$ the smallest natural number $n$ with
the property that the $n$-th up to the $(n+9)$-th digit in the
decimal expansion of $\pi$ form the sequence $0123456789$, and we
define as follows a point $r$ of the reduced continuum: the $n$-th
$\lambda $-interval $\lambda _n$ is a $\lambda ^{(n-1)}$-interval
centered around 0, as long as $n<K_1$; however, for $n\geq K_1$
$\lambda _{n}$ is a $\lambda ^{(n-1)}$-interval centered around
$(-2)^{-K_1}$. The point core of the reduced continuum generated by
$r$ is neither =0, nor $<0$ nor $>0$, as long as the existence of
$K_1$ neither has been proved nor has been proved to be absurd. So
until one of these discoveries has taken place the reduced continuum
is not completely ordered (\cite{Brouwer1991}, p.~31--32).
\end{smquote}

Given an algorithm for $\pi$ all values of this sequence are fixed, so $r$ belongs to the reduced continuum. Note the role of time in this example. Neither $r>0$, nor $r=0$ nor $r<0$ did hold for Brouwer \emph {then and there}. Now we know that $r>0$, see \cite{Borwein1998}. But a similar example is easily constructed.

This example was not new, Brouwer already used it in 1923 (\cite{Brouwer1975}, p.280). But the next example was new, and it is special. Brouwer shows that the natural order $>$ is not an order at all on the full continuum, i.e. he defines a real number $s$  for which neither $s>0$ nor $s<0$ holds, but $s\neq 0$ does hold: 
\begin{smquote}
Therefore we consider a mathematical entity or species $S$, a
property $E$, and we define as follows the point $s$ of the
continuum: the $n$-th $\lambda $-interval $\lambda _n$  is a
$\lambda ^{(n-1)}$-interval centered around $0$, as long as
neither the validity nor the absurdity of $E$ for $S$ is known,
but it is a $\lambda ^{(n)}$-interval centered around $2^{-m}$
($-2^{-m}$), if $n\geq m$ and between the choice of the $(m-1)$-th
and the $m$-th interval a proof of the validity (absurdity) of $E$
for $S$ has been found. The point core belonging to $s$ is $\neq
0$, but as long as neither the absurdity, nor the absurdity of the
absurdity of $E$ for $S$ is known, neither $>0$ nor $<0$. Until
one of these discoveries has taken place, the continuum can not be
ordered ({\cite{Brouwer1991}}, p.~31--32).

\end{smquote}

 We interpreted the sequence $\xi_{\omega}$ in the proof of the NCTH in \cite{Brouwer1927B} as a choice sequence, as is common practice in intuitionistic research. But Brouwer did not mention this special character explicitly in \cite{Brouwer1927B}. That is different in the Berlin Lectures. Because of the distinction reduced continuum versus full continuum, the sequence $s$ must be a choice sequence.  And this time Brouwer has a method to make the application powerful, which is the connection of the values of the terms with the results of working on an untested proposition. Thus Brouwer avoid the contradictorily flavor we noted in the definition of the sequence in the NCTH proof, in always choosing a value from $\xi _0$ but leaving open the possibility of choosing differently. A result on the untested proposition, determining the values of the defined sequence, is in this description  an event arriving from outside.

  Brouwer used the Berlin example many times after this first appearance, for example in his lectures in Vienna in 1928, Geneva in  1934 and in Cambridge after World War II. That Brouwer's definition is a description of a construction is even more clear in his use of the  example in \cite{Brouwer1934}, which is the unpublished text of Brouwer's Geneva Lectures of 1934. Brouwer defines the same sequence $s$ as follows:

 \begin{smquote}
 
 The $n$-th interval $\lambda ^{n-1}$ is of length $2^{-(n-1)}$ and it is centered around the origin. This is how one starts. At the same time one starts to work on a difficult problem, to know whether the property $E$ is true for a space $S$, for example the problem of Fermat. If one finds a solution of this problem between the $(n-1)$-th and the $n$-th choice, then, to start with the $n$-th interval, one will choose that interval in a different way. 
 If the property is true for the space $S$, then for $v \geq n$ the $v$-th interval $\lambda ^{v}$ will  be centered around $2^{-n}$. The next interval will be placed according to this law, in his predecessor and with the same centre.
 When contrary, one would find that the property $E$ is absurd for the space $S$ then the intervals will be centered around $-2^{-n}$.\endnote{Le $n$-ieme intervalle $\lambda^{(n-1)}$ est de longueur $2/(n-1)$ et il est centr\'e a l'origine. C'est ainsi qu'on commense. Mais en m\^{e}me temps on se pose un probl\`{e}me difficile, a savoir si la propriete $E$ est vraie pour une espace $S$, par example le probl\`{e}me de Fermat. L'on trouvait une solution de ce probl\^{e}me entre le $(n-1)$-ieme choix et le $n$-ieme, alors a commencer du $n$-ieme intervalle, on choisirait cet intervalle d'une autre maniere.
Si la propri\'{e}t\'{e} est vraie pour l'espèce $S$, alors le $v$-ieme intervalle sera pour $v\leq n$ l'intervalle $\lambda^{(v)}$ centre au point $2^{-n}$. L'intervalle suivant sera situ\'{e} en vertu de cette loi, dans le pr\'{e}c\`{e}dent et aura le m\^{e}me centre.
Si, au contraire, on trouvait que la propri\'{e}t\'{e} $E$ est absurde pour l'espece $S$, alors centrerait les intervalles autour du point $-2^{-n}$.
Ce point $s$ est defini\'{e} \^{a} une maniere parfaitement correcte. Ce point est different de $0$, puisque si $s$ etais egal a $0$, la possibilite que la suite d'intervalle se continue autour d'un point $2^{-n}$  quelconque, serait exclu. C'est dire que la supposition qu'on un jour la verit\'e de $E$ pour $S$, serait absurde, et que la supposition qu'on trouve un jour l'absurdite pour $S$                  serait aussi absurde. La v\'{e}rit\'{e}é et l'absurdit\'{e} de cette propriet\'{e} seraient absurdes toutes les deux, et cela est impossible. Le point $s$, defini de la sorte, n'est pas egal a $0$, puisqu'il est impossible que la premiere variante de la construction soit exclue en m\^{e}me temps que la seconde.}
\end{smquote}

There can be no doubt that Brouwer above defines a choice sequence, because after the definition Brouwer remarks:

\begin{smquote}
 The point above is not a
sharp point, because the construction is not completely determined,
but depends on the intelligence of the constructor relative to the
posed problem. \endnote {Le point de tout a l'heure n'est pas un point predestin\'{e}, puisqu'il reste dans sa construction quelque chose qui n'est pas entierement determine mais dependre de l'intelligence du constructeur relativement au probleme pose.}
\end{smquote}

Brouwer continues with proving that the reduced continuum is not completely ordered by $<$. In order to do so  he defines the lawlike real $r$, based on the expansion of $\pi$, similarly as in the Berlin Lectures. About the difference between the lawlike $r$ and the choice sequence $s$, he remarks :

\begin{smquote}
	When one hundred different persons are constructing the number
	$r$, one is always certain that any interval chosen by one of
	these persons is always covered, at least partly, by every
	interval chosen by one of the others. That is different for $s$. If I would give the definition of $s$ to one hundred persons, who are all going to work in a different room, it is possible that one of these one hundred persons at one time will
choose an interval not covered by an interval chosen by one of the
others.\endnote{`Si cent personnes differentes  s'occupent de la construction du nombre $r$, on est toujours sur que chaque intervalle choisi par une de ces personnes, sera couvert, du moins a partie, par chaque intervalle choisi par une autre de ces personnes. Il n'est pas de m\^{e}me pour $s$. Si je donne la definition de $s$ a cent personne differentes, qui travaillent dans un local different, il se peut que l'une de ces cent personnes choisisse une fois un intervalle ne couvrait pas l'intervalle choisi par une autre de ces personnes.'}
\end{smquote}

\bigskip

We think the quotes from the Geneva Lectures are important for understanding the concept of choice sequence. Nowhere in his work Brouwer is so clear about this concept and as far as we know it is the only place in the work of Brouwer where he explicitly mentions a difference between a lawlike sequence and a particular choice sequence. In constructing a lawlike sequence different constructors reach the same values for the terms. In case of a choice sequence this does not have to be the case.

Our conclusion is: a choice sequence is given by a description of its construction. The description does not determine the values of the terms completely. In the example Brouwer uses these values are made to depend on the future mathematical experience of the constructor. The reasoning is done on the basis of the description only, before the construction has started.

A consequence is that for the choice sequence $s$, used in the Berlin and Geneva Lectures and quoted above, and for some natural number $n_0$, we do not have $s_{n_0}=0 ~ \vee ~ s_{n_0} \neq 0$: we reason before the construction has started and we do not know whether and when we have a result on working on the specific problem. As far as we know $a_n=0 ~\vee ~ a_n \neq 0$ is  valid for all $n$ in all existing formalization's of choice sequences, which makes them unfit for a description of Brouwer's choice sequences. An argument for $s_{n_0}=0 \vee s_{n_0} \neq 0$ is that we can go to stage $n_0$. But different persons may reach $n_0$ with different values. What holds for $s$ must hold for all the possible continuations of the use of its definition.

In the next section we shall present principles we derive from the conception of  choice sequence above in order to analyse Brouwer´s CS arguments.

\section{The theory of individual choice sequences TICS}

In his first publication after the World War 2, \cite{Brouwer1948A}, Brouwer uses the Berlin example again. Brouwer had switched from German to English. Whereas  Brouwer had used `we' above, uses `I' in a Vienna example below and  `constructor' in the Geneva lectures to denote the constructor of the sequence generating $s$, he now uses the expression `creating subject' for that purpose. Although Brouwer remarked that he had used this example in his lectures since 1927, the applied method has mistakenly been supposed to be radically new in 1948, see e.g. the quote of Heyting in the introduction. 

The  expression `the creating subject', less appropriately renamed `creative subject',  was interpreted as `the idealized mathematician', the IM, his whole mathematical activity covered by an $\omega$ sequence of discrete stages. The results on these stages determine the values of the defined sequence. Kreisel formalized the method according to this conception, see \cite{Kreisel1967}, and Troelstra's elaborated this formalization in \cite{Troelstra1969}, setting the standard.

The basic term of the TCS is $\Box _n \varphi$, were $\varphi$ is a mathematical assertion, expressing that at stage $n$ the IM has evidence for $\varphi$. We shall use the language of the TCS to formulate the principles for  handling Brouwer's argument derived from the conception of the previous section, with a different interpretation of the basic term. Thus, we have the possibility to compare reconstructions and point out differences. Furthermore, this formalization leads to new insights.

\smallskip

The crucial difference between our interpretation and the TCS is that in the TCS  the stages embody the whole mathematical activity of the IM. In our interpretation they cover only the future. All stages lie before us at the moment of reasoning about the choice sequence. We define for a formula $\varphi$ 

\bigskip
	
	$\Box_n \varphi$ 

 \bigskip

as `at the $n$-th stage from now $\varphi$ will hold'.  $\Box _n \varphi$ may hold because $\varphi$ holds now, but a proof of $\Box_n \varphi$ may depend also on information coming available at or before stage $n$. Information does not get lost. Once a proof is available it will remain so. The stages may be connected to the decisions taken about the values of a particular choice sequence on the basis of the work on an untested proposition.

	The time it takes for the proof to be carried out does not play a role. For example, for the lawlike real $r$ in the previous section and for some $n_0$ we have $r_{n_0}=0~ \vee~ r_{n_0} \neq 0$, because all information to calculate $r_{n_0}$ is available here and now. Like Brouwer, we completely disregard the time or complexity of the necessary computations.
 
 For the choice sequence $s$ in the previous section we do not have $s_{n_0}=0 ~ \vee ~ s_{n_0} \neq 0$, because we do not now whether we shall have reached a result concerning the untested proposition before stage $n_0$, if we let the stages coincide with the choices. But at stage $n_0$ we know whether we have reached a result on stage $n_0$, so we have   $\Box _{n_0} (s_{n_0}=0 ~ \vee ~ s_{n_0} \neq 0)$.
	
	Under this interpretation of $\Box _n$ we adopt for a formula $\varphi$ the following principles:
	
	\bigskip

IC1 ~~ $\Box_n \varphi \rightarrow  \Box_{n+m} \varphi$,

\smallskip

Expressing that once $\varphi$ has been proved, it will remain valid in the future.

\smallskip

IC2 ~~ $\neg \varphi \rightarrow \neg \exists n \Box_n \varphi$.

\smallskip

Expressing that if the assumption of a proof of $\varphi$ is contradictory, the assumption of a future proof is also contradictory.

\smallskip

IC3 ~~  $\varphi \rightarrow \exists n \Box_n\varphi$,

\smallskip

For, if we have a proof of $\varphi$ now, IC3 holds by definition. And when we will acquire a proof $\varphi$ in the  future, it will be in a certain stage, say $n_0$, which gives us immediately $\Box_{n_0} \varphi$. But Brouwer does not fully use IC3. His main device is, as we shall see in the examples to come, the (weaker) contraposition of IC3, i.e.

\smallskip

MD  $ \neg \exists n \Box_n\varphi \rightarrow \neg \varphi$.

\smallskip

Brouwer's main device MD is the only principle for which Brouwer explicitly gives a motivation. In one of the rare occasions that Brouwer gives a detailed proof, see \cite{Brouwer1948A}, he reasons for an untested assertion $\alpha$ that if it is certain certain that $\alpha$ could never be proved to be true, the the absurdity of $\alpha$ would be known. It might well be that the `enormous discovery` from the quote in section 6 refers to this principle.

\smallskip

The following two principles form with IC1 and IC3 the standard theory TCS.

\bigskip
  
CS4 ~~ $\Box _n \varphi \vee \neg \Box _{n}  \varphi $ for each $n$, 

\smallskip

CS5 ~~ $ \exists n \Box_n \varphi \rightarrow \varphi$.

\bigskip

In our interpretation CS4 expresses that for every $n$ and $\varphi$  we now have a proof that $\varphi$ will hold at stage $n$, or that the supposition that $\varphi$ will hold at stage $n$ is contradictory, which is of course not valid; we reject CS4.

Above we argued that for any $n_0$ we have  $\Box_{n_0} (s_{n_0}=0 \vee s_{n_0} \neq 0)$ but not $(s_{n_0}=0 \vee s_{n_0} \neq 0)$, which is a counterexample to CS5.  In this counterexample we use a term of a sequence which is  yet undetermined, i.e. $\varphi$ contains a choice parameter. If $\varphi$ does not contain choice parameters we see no objection against CS5, and we will accept CS5 in this restricted form, analogous to the acceptance of $(r_n=0 \vee r_n \neq 0)$ for the special case that $r$ is lawlike. But we reject CS5 for the general case. Consequently, for general $\varphi$ we do not have $\exists n \Box_n (\neg \varphi \vee \neg \neg \varphi) \rightarrow (\neg \varphi \vee \neg \neg \varphi)$. This corresponds to a distinction Brouwer actually makes in his Post-War papers. We shall have the occasion later on to discuss \textit{}this further. 

 \smallskip
 
By defining
\smallskip

$a_n=0$ if $\Box _n \varphi$, and $a_n=1$ if $\neg \Box _n \varphi$

\smallskip

we can deduce in the TCS for any mathematical assertion $\varphi$ the following schema, called Kripke's Schema

 \bigskip

 {\bf KS~~~} $\exists a(\forall n( a_n=0
\vee a_n=1)~ \wedge ~ \exists n(a_n\neq 0 \leftrightarrow \varphi))$.

\bigskip

 KS is  is sufficient for most of Brouwer's CS-arguments. It has the advantage that it does not contain the problematic notion of IM, and can therefore be used for Brouwer's CS arguments instead of the TCS. Since for the derivation of KS the principles CS4 and CS5 are used, which we reject, we have no reason to accept KS.

\bigskip

In our \cite{Niekus1987} we presented a Kripke model in which the principles IC1, IC2 and IC3 all are valid, and CS4 and CS5 are not. We claim that the model is sufficiently close to Brouwer´s thoughts to discuss his arguments and to clarify them, but there is no claim that the model exactly represents his ideas.  The formulas valid in these models form the Theory of Individual Choice Sequences TICS. We propose to replace the standard theory TCS by our TICS to analyse  CS arguments. 
\bigskip

In the coming sections we shall consecutively discuss CS arguments from the Vienna Lecture in 1928 (\cite{Brouwer1930A}, Amsterdam 1946 (\cite{Brouwer1948C}) and his Canadian Lecturers from 1953 (\cite{Brouwer1954A}). In these lectures Brouwer presents new results on the continuum based on particular choice sequences. The terms of these sequences are chosen from a converging sequence of real numbers and its limiting number, depending on the results of the constructor working on an unsolved problem. Only in the text of the final lecture will he give proofs for his results. These proofs are important for our position.
 
 We reconstruct his intended argument with IC2, IC3 and the restricted form of CS5. It will be clear that Brouwer is very aware of the different possibilities inherent in CS5 by the assumptions he makes. We shall show that the acceptance of TC5 for the special case, and its rejection of the general case is fully supported by Brouwer's practice. We will conclude that Brouwer's CS arguments give no basis to accept KS.

\section{ Vienna 1928 }

A year after Berlin,  Brouwer lectures in Vienna. In the second of his two lectures the new method introduced a year before in Berlin comes to full maturity (\cite{Brouwer1930A} ). Brouwer  examines the continuum regarding seven properties, all valid classically but not intuitionisticly (\cite{Brouwer1930A}). Each time he carefully distinguishes whether the result counts only for the full continuum, or for the reduced continuum as well, i.e. whether a choice sequence is used in the counterexample or a lawlike one. 

Instead of using an arbitrary untested `property $E$ for a species $S$´   as he did in the Berlin Lectures, this time Brouwer uses the newly introduced notion of a fleeing property, a generalization of the technique with the decimal expansion of $\pi$.

\bigskip

\emph{A fleeing property for natural numbers} satisfies the following conditions: it is decidable for each $n$, no natural number is known for which the property holds, and the assumption of the existence of a natural number possessing the property is not known to be contradictory. The critical number $\lambda _f$ (Lösungszahl) of a fleeing property $f$ is the smallest natural number possessing the property $f$.  
\bigskip

Brouwer's standard example of a fleeing property is being the first digit in the decimal expansion of $\pi$ in a row of $1$, $2$, ... $9$ (which is for us nowadays not fleeing anymore). His  definition of a fleeing property does not exclude that the contradictoriness of the existence of a natural number possessing the property, itself is contradictory, i.e. testedness. But if it would be contradictory, in some of his examples one of the possibilities in the definition would be excluded, so one may fairly assume that Brouwer presupposed untestedness. In his Post-War papers Brouwer would add this condition in using fleeing properties when needed. He even attempts to use a fleeing property for which the absurdity of the absurdity of its existence has been established, see section 13.

 In the second Vienna Lecture Brouwer starts with the same two examples as in the Berlin Lecture. Then he introduces the virtual order, see section 4.  We repeat: the relation $x\prec y$ defined by $x \prec y \leftrightarrow (x \neq y \wedge \neg x > y)$ satisfies the properties of a virtual order. With respect to this virtual order Brouwer proves that the continuum is not dense in itself.

\smallskip

A species $S$ is \emph{dense in itself} (in sich dicht) if every point of $S$ is a limiting point.
\smallskip

An element $a$ of a species $S$ is a \emph{limiting point} (Grenzelement) if for a sequence $a_1, a_2, ...$ with $a_1 \prec a_2 \prec .. \prec a$  for every $b  \prec a$ an $a_n$ exists such that $a_n \succ b $, or analogously for a decreasing sequence.

\smallskip

Brouwer's proof  that the full continuum is not dense in itself: 

\bigskip

Let $a_1,  a_2,...$ be such that $a_1 \prec  a_2  \prec...1/2$. Let $f$ be a fleeing property and let $\lambda_f$ be its critical number. We define a sequence  $e_1,~e_2, ...$ as follows:
As long as we did not find the critical number $\lambda_f$, nor the contradictory of its existence has become known,   we choose $e_n=a_n$, but as soon as between the choice of $e_v $  an $e_{v+1}$ one of these two events happens,  we choose $e_{v+n}=a_v$ for all n.
The sequence $e_1,~e_2,...$ defines an element of the continuum, say $e$. Then $e\neq 1/2$, but there is no $m$ such that $a_m \succ e$.

\bigskip

Brouwer gives no further argument for the results. We shall give a reconstruction in the language of the previous section. We express  $\exists n ( \lambda _f=n)$  by $\alpha$. The crucial step is Brouwer's main device MD.

\bigskip

From the definition follows $\exists n \Box _n (\alpha \vee \neg \alpha) \leftrightarrow \exists v (e=a_v)$.

\smallskip

Suppose $e = 1/2$. 

\smallskip

Because $e=1/2 \rightarrow \neg \exists v ( e=a_v)$,

\smallskip

we would have $\neg \exists n \Box _n(\alpha \vee \neg \alpha)$.
\smallskip 

Then according to MD we would have $\neg( \alpha\vee \neg \alpha)$, a contradiction.

\smallskip
So $e \neq 1/2$.
\bigskip

Suppose $\exists v(e \prec a_v )$.
\smallskip

So $\exists n (a_n=e)$ holds, and we would have $\exists n \Box_n (\neg \alpha \vee \alpha)$. 

\smallskip
Since $\alpha$ expresses the existence of a fleeing property, $\alpha$ does not contain a choice parameter. So, as we argued in the previous section, we may apply TC5, from which it follows that $\alpha \vee \neg \alpha$ has to hold, which it does not by assumption.

\smallskip

So $\exists v(e \prec a_v )$ does not hold.

\smallskip

We shall see in the next section that Brouwer refines his definition to avoid the use of TC5 for an arbitrary untested proposition.

\bigskip

The resemblance between the sequence used in the proof of the NCTH and the one above is striking. In both cases the values of the defined sequence are chosen from a convergent sequence of reals and its limiting number. But in contrast with the proof of the NCTH, the Vienna proof, published three years later, did not receive any attention until very recently (\cite{Fletcher2020}). 

It is hard to understand that \cite{Brouwer1930A}, the text of the second Vienna Lecture, has never been mentioned in twentieth century research concerning Brouwer's choice sequences. A fleeing property is a generalization of a technique used by Brouwer to define lawlike sequences. Its  use in defining a choice sequence may in a superficial reading have made it difficult to recognize its special character. But given the use of the distinction full continuum versus reduced continuum, and the explicit statements of results as applying to the full continuum only or for both, there can be no misunderstanding that Brouwer uses individual choice sequences extensively.  

After each proof with a choice sequence Brouwer adds a puzzling remark, namely that the same holds for the reduced continuum. Then why did he not use a lawlike sequence, as for example in the proof that the reduced continuum is not dense in itself, see above? We shall return to this question in section 13.

Finally, it is remarkable that Brouwer after this first publication of the breakthrough on particular choice sequences, and his extensive exploration of the distinction reduced versus full continuum, at that moment unique in his work, did not publish anything seriously for more than fifteen years. We do not think that this inactivity had its cause in Brouwer's dissatisfaction with his proof of the fan theorem which is the view of van Stigt we mentioned in section 6. More likely it was caused by a conflict between Brouwer and Hilbert, shortly after the Vienna lectures. It made Brouwer lose his position in the editorial board of the \emph{Mathematische Annalen}; it touched him deeply, see \cite{Dalen2005}, pages 599 and next.

When Brouwer started to publish again after the Second World War, he continued, in so far as it concerned choice sequences, exactly where he had left off in 1930.

 \section{The Post War papers - a refinement}

In his first Post War  publication, \cite{Brouwer1948A}, Brouwer treats again the first Berlin choice sequence,  this time with a rare detailed proof. We treated this example in our former papers. Here our focus will be on \cite{Brouwer1948C} and \cite{Brouwer1954A}. In these papers   Brouwer gives his most sophisticated counterexamples with individual choice sequences, using a newly introduced concept: the notion of a drift. Our reconstruction enables a comparison with reconstructions in the recent \cite{Atten2018} and \cite{Fletcher2020}. But above all, \cite{Brouwer1954A} is the other rare occasion in which Brouwer gives proofs. One of these proofs is used to support the standard theory TCS. We claim, to the contrary, that our analysis provides convincing evidence against the TCS.  

\smallskip

   In the definition of a choice sequence in \cite{Brouwer1948C} and \cite{Brouwer1954A} Brouwer adds to the condition  `$\alpha$ is not tested' also the condition `$\alpha$ is not recognized as testable'. Clearly  Brouwer introduces this new condition to express a real distinction and in our reconstruction we see that Brouwer has  a  good reason for its introduction. As we will see it is here that the distinction between the restricted and unrestricted versions of TC5 comes out very clearly.
   
   The proper candidate to express `recognized as testable' in our formalism must be $\exists n \Box_n (\neg \alpha \vee \neg \neg \alpha)$. In section 8 we argued that $s_{n_0}=0 \vee s_{n_0} \neq 0$ for the choice sequence $s$ does not hold because we do not know now whether we will have at stage $n_0$ established  a result on the posed problem.  But at stage $n_0$  we will know, so we do have $\Box_{n_0}(s_{n_0}=0 \vee s_{n_0} \neq 0) $  which is a counterexample to CS5, i.e. we reject ~~   [CS5 $\exists n \Box_n \alpha \rightarrow \alpha$] ~~ for the general case. 
   
   We see no reason to reject CS5 if $\alpha$ does not contain choice parameters. In the example of the previous section we saw that Brouwer concludes from $\exists n \Box_n (\alpha \vee \neg \alpha)$ that $\alpha$ is tested, which is correct, since $\alpha$ is the assertion $\exists n (\lambda_f = n)$ for a fleeing property $f$, which does not contain a choice parameter. In the  examples we will discuss in sections 11 and 12 Brouwer cannot conclude  $\alpha$ is tested from $\exists n \Box_n (\neg \alpha \vee \neg \neg \alpha)$ since $\alpha$ is arbitrary,, so he adds `$\exists n \Box_n (\neg \alpha \vee \neg \neg \alpha)$ does not hold' as a condition on the $\alpha$ used.

In a note found with the manuscript of \cite{Brouwer1949A} Brouwer gives a notation for this  distinction with the introduction of  $|a$ and $a|$. The former expresses `$a$ is tested from now on',  and the later expresses `$a$ is tested from a certain moment in the future', which are of course nothing but our $(\neg a \vee \neg \neg a)$ and $\exists n  \Box_n (\neg a \vee \neg \neg a)$. Brouwer argues in the note that the former is contradictory and the latter is not, so he  rejects CS5 in a very strong way. The $\alpha$ used in this argument is the statement that a sequence, now completely free, defines a rational real, so it clearly contains a choice parameter.

At the end of this note Brouwer expresses doubts about introducing $|a$ as a mathematical notion, without further argument.  Brouwer's doubts and change of position is characteristic for his struggle with the notion of choice sequence. But since Brouwer's $|a$ is just our $\neg a \vee \neg \neg a$  acceptance of CS5 is of course out of the question.

\section{The notion of a drift}

In his \cite{Brouwer1927B} proof of the NCTH Brouwer defines a sequence whose terms are chosen from a convergent sequence of reals and its limiting number. In his Berlin Lectures he makes this method work by connecting the choices with the results of the constructor working on an unsolved problem. In his Post War papers he refines this technique further with the introduction of the notion of a drift. This notion enables Brouwer to express subtle distinctions of intuitionistic logic which do not exist classically.

\smallskip

A \emph{drift} $\gamma$ is the union of an infinitely preceding  convergent sequence of real numbers $c_1,~c_2, ...$, the \emph{counting-numbers} of the drift, with the limiting number $c$, named the \emph{kernel} of the drift; all counting numbers lie apart from each other and from the kernel. If $c_v > c$ for each $v$ the drift is \emph{right-winged}, if $c_v < c$ for every $v$, the drift is \emph{left-winged}. If $c_1,~c_2,...$ is the union of a sequence of left counting-numbers $l_1, l_2, ...$ such that $l_v< c$ for each $v$, and a sequence of right counting-numbers $r_1,~r_2,...$ such that $r_v > c$ for each $v$, the drift is \emph{two-winged}.

\bigskip

In \cite{Brouwer1948C} Brouwer gives the following application. Let $\alpha$ be a mathematical assertion neither tested nor recognized as testable. In connection with $\alpha$ and a drift $\gamma$ we can define an infinitely proceeding sequence $d_1,~d_2,...$ of real numbers  as follows: 

As long as during the choice of the $d_n$ we have neither experienced the truth nor the absurdity of $\alpha$, each $d_n$ is chosen equal to $c$.

But as soon as, between the choice of $d_{r-1}$ and $d_r$ we experience the truth or the absurdity of $\alpha$,  $d_r$ as well as $d_{r+v}$ for all $v$, is chosen to be $c_r$. This sequence $d_1,~d_2,...$ converges to a real number $d_{\gamma , \alpha}$, which is called the \emph{direct checking-number of $\gamma$ through $\alpha$}.

\bigskip
Brouwer states: Let $\gamma$ be a right-winged drift whose counting-numbers are rational. Then the rationality of $d_{\gamma , \alpha}$ is testable but not decidable, and its non-contradictorily is not equivalent to its truth. Furthermore we have $d_{\gamma , \alpha} \succ c$, but not $d_{\gamma , \alpha} > c $.

Brouwer does not give a proof. 

Our reconstruction, with $R(x)$ expressing `x is rational':

\bigskip

 According to the definition of $d_{\gamma , \alpha}$ we have
 
 \smallskip
 
   $R(d_{\gamma , \alpha}) \leftrightarrow  \exists v (c_v=d_{\gamma , \alpha})$ and $ \exists v (c_v=d_{\gamma , \alpha}) \leftrightarrow \exists n \Box_n(\alpha \vee \neg \alpha) $.
   
   \smallskip
   
  Since $d_{\gamma , \alpha}=c \rightarrow \neg \Box(  \alpha \vee \neg \alpha)$  and $\neg \Box ( \alpha \vee \alpha \neg \alpha) \rightarrow \neg (\alpha \vee \neg \alpha)$, 
  
  by MD we have $d_{\gamma, \alpha} \neq c$.
\textit{}
   
\smallskip

   Now suppose $\neg R(d_{\gamma , \alpha})$. 
   
  \smallskip

   Then $\neg \exists v (c_v=d_{\gamma , \alpha})$,
   
   \smallskip
   
   so $\neg \exists n \Box_n(\alpha \vee \neg \alpha)$, and consequently (MD),  $\neg (\alpha \vee \neg \alpha)$, a contradiction.
   
   \smallskip
   
   So $\neg \neg R(d_{\gamma , \alpha})$, and a fortiory, $R(d_{\gamma , \alpha})$ is tested.
   
   \smallskip
   
   If $R(d_{\gamma , \alpha})$ then $\exists n \Box_n (\alpha \vee \neg \alpha)$, and $\alpha$ would be recognized as testable, which 
   
   it is not. 
   So $R(d_{\gamma , \alpha})$ does not hold, and consequently also 
   
   $\neg \neg R(d_{\gamma , \alpha}) \leftrightarrow R(d_{\gamma , \alpha}) $ is not valid.
   
   \bigskip
   
   By definition  $d_{\gamma , \alpha} \succ c \leftrightarrow (d_{\gamma , \alpha} \neq c \wedge \neg d_{\gamma , \alpha}<c $).
   
   \smallskip
   
   Because $\gamma$ is a right-winged drift $\neg d_{\gamma, \alpha} < c$ and $d_{\gamma , \alpha} \neq c$,
   we have $d_{\gamma , \alpha} \succ c$.
   
   \smallskip   
   
   If $d_{\gamma , \alpha}>c$ we would have $\exists n \Box_n(\alpha \vee \neg \alpha)$, and $\alpha$ would be recognized as 
   
   testable, which it is not, so $d_{\gamma , \alpha}>0$ does not hold. 

   \bigskip

Brouwer gives another application of a drift in \cite{Brouwer1948C}. First he defines:

\smallskip

Again, in connection with an assertion $\alpha$ neither tested nor recognized as testable  and a two-winged drift $\gamma$ we can construct an infinitely proceeding sequence $e_1,~e_2,...$  of real numbers  as follows: 

As long as during the choice of the $e_n$ we have  neither proved $\alpha$ or $\neg \alpha$ , each $e_n$ is chosen equal to $c$. 

But a soon as between the choice of $e_{v-1}$  and $e_v$  we have obtained proof of $\alpha$, then for each $w$, $e_{v+w}$ is chosen to be $r_v$. And as soon as between the choice of $e_s$  and $e_{s+1}$  we proved the absurdity of $\alpha$,  for all $v$, $e_{s+v}$ is chosen to be $l_s$. This sequence $e_1,~e_2,...$ converges to a real number $e_{\gamma , \alpha}$, the \emph{oscillatory checking-number of $\gamma$  through $\alpha$}.

\smallskip

Then he states:

\bigskip

Let $\gamma$  be a two-sided drift whose right counting-numbers $r_1,~r_2,...$ are all rational and whose left counting-numbers $l_1,~l_2,...$ are all irrational. Then the rationality of $e_{\gamma , \alpha}$ is neither decidable nor testable, nor is its non-contradictory equivalent to its truth. Finally, neither $e_{\gamma , \alpha} \succeq c$ nor $c \succeq e_{\gamma , \alpha}$ does hold.

\smallskip
We repeat that Brouwer introduced $x \prec y$ by $x \neq y \wedge \neg x > y$. Conform this definition he introduces in this paper $x \succeq y$ by $\neg x<y$.

\bigskip
Brouwer does not give a proof. Our reconstruction: 

 \smallskip
 First, $c \neq e_{\gamma, \alpha}$ because

 \smallskip

 $e_{\gamma, \alpha}=c \to ( \neg \exists n \Box_n \alpha \wedge \neg \exists n \Box_n \neg \alpha)$, and 

 \smallskip
 
 $( \neg \exists n \Box_n \alpha \wedge \neg \exists n \Box_n \neg \alpha) \to (\neg \alpha \wedge \neg \neg \alpha)$,  a contradiction.

\smallskip 

So $R(e_{\gamma, \lambda}) \leftrightarrow \exists v (e_{\gamma , \alpha}= r_v) \leftrightarrow \exists n \Box_n \alpha$,

\smallskip

 $ \neg R(e_{\gamma, \alpha)}\leftrightarrow  \neg\exists n \Box_n \alpha \leftrightarrow \neg \alpha$. 

\smallskip

Analougosly $\neg R(e_{\gamma, \alpha}) \leftrightarrow \exists v (e_{\gamma , \alpha}= l_v) \leftrightarrow \exists n \Box_n \neg \alpha$,

\smallskip

 $\neg\neg R(e_{\gamma, \lambda}) \leftrightarrow \neg \neg \alpha$.

\smallskip
Since $\alpha$ is not tested, neither is the rationality $R(e_{\gamma , \alpha})$,

\smallskip

and consequently, neither is the rationality decidable.

\bigskip

$R(e_{\gamma , \alpha})$ and $\neg \neg R(e_{\gamma , \alpha})$  are not equivalent: 

\smallskip 

$R(e_{\gamma , \alpha}) \leftrightarrow \exists v e_{\gamma, \alpha}=r_v$     and $\exists v (e_{\gamma, \alpha}=r_v) \leftrightarrow e_{\gamma , \alpha}>c $;

\smallskip

$\neg \neg R(e_{\gamma , \alpha}) \leftrightarrow \neg \exists v ( e_{\gamma , \alpha}=l_v)$ and $\exists v (e_{\gamma , \alpha} l_v) \leftrightarrow \neg e_{\gamma ,\alpha} < c $.

For $R(e_{\gamma , \alpha})$ to hold, a positive distance between $e_{\gamma , \alpha}$ and $c$ must exist, 

for $\neg \neg R(e_{\gamma , \alpha}) $ to hold, this does not have to be the case.

 \bigskip
 Finally, because of $\neg e_{\gamma , \alpha} < c \rightarrow \neg \exists n \Box _n \neg \alpha$ and $\neg \exists n \Box _n \neg \alpha \rightarrow \neg \neg \alpha$,
 
 \smallskip
  and because of $\neg e_{\gamma , \alpha} > c \rightarrow \neg \exists n \Box _n  \alpha$ and $\neg \exists n \Box _n \neg \alpha \rightarrow  \neg \alpha$,
 
 \smallskip
 
 neither $e_{\gamma , \alpha} \succeq c$ nor $c \succeq e_{\gamma , \alpha}$ can hold, since $\alpha$ would then be tested.

 \section{Kripke's Schema}

The notion of a drift appears again in \cite{Brouwer1954A}. After the direct and the oscillating checking number, Brouwer defines a third checking number. For a drift $\gamma$ with counting numbers $c_1,~c_2, ...$  and kernel $c$, and for a mathematical assertion $\alpha$ which is neither tested nor recognized as testable, Brouwer defines the conditional checking number $f_{\gamma , \alpha}$ as follows.

\smallskip
We can generate an infinitely proceeding sequence $f_1,~f_2,...$ by choosing $f_n=c$ as long as we did not find a proof for $\alpha$. But as soon as, between the choice of $f_r$ and $f_{r+1}$ we have found a proof for $\alpha$,  we choose $f_{r+v}= c_r$ for all $v$. The sequence $f_1,~f_2,...$ converges to a limiting number, the \emph{conditional checking number} $f_{\gamma , \alpha}$.

\smallskip

So the difference between $f_{\gamma , \alpha}$ and the direct checking number $d_{\gamma , \alpha}$ is that the values of the generating sequence of $d_{\gamma , \alpha}$ change when we experience the truth or the absurdity of $\alpha$, while in case of $f_{\gamma , \alpha}$ the values change only if we experience the truth of $\alpha$. 

Brouwer examines the consequences of this difference for a drift $\gamma$ with rational counting numbers and an irrational kernel, and an assertion $\alpha$ neither tested nor recognized as testable. For this drift Brouwer states that the rationality of the direct checking number $d_{\gamma,\alpha}$  is testable, but not decidable (judgeable). Brouwer's argument that it is testable :

\smallskip

`Then the assertion of the  rationality of a $d_{\gamma , \alpha}$ is not judgeable, but it is testable, because the assertion of irrationality of $d_{\gamma , \alpha}$ would entail the simultaneous contradictory of the truth and the absurdity of $\alpha$,  which is an absurdity.' 

\smallskip
 
 If again $R(x)$ expresses `x is rational' we can reconstruct Brouwer's reasoning in deriving a contradiction as follows:
 
 \smallskip

 $R(d_{\gamma , \alpha}) \leftrightarrow \exists n \Box _n(\alpha \vee \neg \alpha)$ by the definition of $d_ \gamma , \alpha$.

 \smallskip
 
Suppose $\neg R(d_{\gamma , \alpha})$. 
 
 \smallskip
 
Because of $\neg R(d_{\gamma , \alpha}) \rightarrow \neg \exists n \Box _n(\alpha \vee \neg \alpha)$, 
\smallskip

  and $ \neg \exists n \Box_n (\alpha \vee \neg \alpha) \rightarrow \neg (\alpha \vee \neg \alpha)$,    Brouwer's main device MD,

  \smallskip  $\neg(\alpha \vee  \neg \alpha)$ were to hold, which is a contradiction.

  \smallskip
  So $\neg \neg  R(d_{\gamma , \alpha})$.

\smallskip

$R(d_{\gamma ,\alpha})$ can not hold because $\alpha$ would have been recognized as testable,

\smallskip

so $\neg \neg R(d_{\gamma ,\alpha}) \rightarrow R(d_{\gamma ,\alpha})$ does not hold either. \\

For the conditional checking number $f_{\gamma ,\alpha}$ Brouwer states that neither $R(f_{\gamma ,\alpha})$, nor $\neg \neg R(f_{\gamma ,\alpha})$, nor $\neg \neg R(f_{\gamma ,\alpha})\leftrightarrow R(f_{\gamma ,\alpha})$ does hold. He starts his proof with  

\smallskip
`On the other hand, truth of $\alpha$ and the rationality of $f_{\gamma , \alpha}$ are equivalent.'
\smallskip

This statement has been interpreted, first by John Myhill (\cite{Myhill1967} p.295), as Brouwer is using an instance of Kripke's Schema KS. We reject, as we explained before, KS because for its derivation the following principles are used: 

 \bigskip{}

CS4~~~  $\Box _n \varphi \vee \neg \Box _{n}  \varphi $ (decidability), and

\smallskip
CS5~~~  $ \exists n \Box_n\varphi \rightarrow \varphi$.

\bigskip{}

In section 7 we showed that neither principle is valid in our interpretation of $\Box_n \varphi$. We also showed that our rejection of CS4 is conform Brouwer's practise. We claim that the same holds for CS5. 

Concerning CS5, we saw a use of his principle in the example of the Vienna lecture in section 8. But the $\alpha $ used did not contain choice parameters, in which case the principle seems acceptable. But for the general case Brouwer evidently introduces the distinction `not recognized as testable' versus `not tested' to avoid its use. In a note which we showed in the previous section Brouwer even gives a special notation for this distinction, the only one in his use of particular choice sequences.  So what is going on with Brouwer's statement that truth of $\alpha$ and the rationality of $f_{\gamma , \alpha}$ are equivalent? A closer look at the complete proof will give the answer. It runs as follows:

\begin{smquote}

On the other hand, truth of $\alpha$ and the rationality of $f_{\gamma , \alpha}$ are equivalent. So the rationality of $f_{\gamma , \alpha}$ is neither decidable nor testable. For, non-contradictory of rationality of $f_{\gamma , \alpha}$ would entail non-contradictory of $\alpha$, i.e. testability of $\alpha$, which was presupposed not to exist. Furthermore, if some day $\alpha$  would prove to be non-contradictory without being true, rationality of $f_{\gamma , \alpha}$ likewise would be non-contradictory without being true. So for rationality of $f_{\gamma , \alpha}$, just as for $\alpha$, non-contradictory would not be equivalent to truth (\cite{Brouwer1954A}, p.4).

\end{smquote}
 
 So Brouwer reasons:
 
 \smallskip
 
 (1) $R(f_{\gamma , \alpha}) \leftrightarrow \alpha$. 
 
 \smallskip
 
 (2) So $R(f_{\gamma , \alpha})$  is neither decidable nor testable. 
 
 \smallskip
 (3) For, $\neg \neg R(f_{\gamma , \alpha}) \rightarrow \neg \neg \alpha$ and $\alpha$ was supposed to be not tested.

  \smallskip
 
 (4) Furthermore if $\exists n \Box_n \neg \neg \alpha$ and not   $\exists n \Box_n \alpha$, 
 
 then also $\neg \neg R(f_{\gamma , \alpha})$ and not $R(f_{\gamma , \alpha})$. 
  
  \smallskip
  (5) So not $\neg \neg R(f_{\gamma , \alpha}) \rightarrow R(f_{\gamma , \alpha})$, just as for $\alpha$. \\

\smallskip

According to us from the definition one may only conclude to   $R(f_{\gamma , \alpha}) \leftrightarrow \exists n \Box _n \alpha$. But Brouwer started with stating $R(f_{\gamma , \alpha}) \leftrightarrow \alpha$, which seems to imply an instance of  $ \exists n \Box _n \alpha \rightarrow \alpha$. But then why Brouwer did not state immediately that $\neg \neg R(x)$, $R(x)$ and  $\neg \neg R(x) \rightarrow R(x)$ all three do not hold? They follow immediately from (1), since from (2) and (5) one may conclude that Brouwer assumes that $\alpha$, $\neg \neg \alpha$ and $\neg \neg \alpha \rightarrow \alpha$ do not hold. And furthermore, then line (3) and (5) are superfluous, and  (4) is incomprehensible.

\smallskip

A closer look at Brouwer's proof shows that though he starts with 

$R(f_{\gamma , \alpha}) \leftrightarrow \alpha$, he reasons as if he had started with 

\smallskip

(1a) $ R(f_{\gamma , \alpha}) \leftrightarrow \exists n \Box _n \alpha$.

\smallskip
Now the proof runs as follows:
\bigskip

$\neg \neg R(f_{\gamma , \alpha})$ does not hold because of

\smallskip

$\neg \neg R(f_{\gamma , \alpha}) \leftrightarrow \neg \neg \exists n \Box_n \alpha$
and $\neg \neg \exists n \Box _n \alpha \leftrightarrow \neg \neg \alpha$ (twice MD),

\smallskip

it would follow that $\neg \neg \alpha$ were to hold, which it does not, 

\smallskip

so $\neg \neg R(f_{\gamma, \alpha})$ does not hold.

\bigskip

Further, suppose $\exists n \Box_n \neg \neg \alpha$ were to hold,  and $\exists n \Box_n \alpha$ was not.

\smallskip

If $\neg R(x)$, then because

\smallskip

$ \neg R(f_{\gamma , \alpha}) \leftrightarrow  \neg \exists n \Box_n \alpha$
and
$ \neg \exists n \Box _n \alpha \leftrightarrow  \neg \alpha$,

\smallskip

$\neg \alpha$ would hold, contradicting $\exists n \Box_n \neg \neg \alpha$,

\smallskip

so $\neg \neg R(f_{\gamma , \alpha})$.

\smallskip
  Since  $ R(f_{\gamma , \alpha}) \leftrightarrow \exists n \Box _n \alpha$, $R(f_{\gamma, \alpha})$ does not hold, so neither does
  
  \smallskip
  
  $\neg \neg R(f_{\gamma , \alpha}) \rightarrow R(f_{\gamma , \alpha})$.

\bigskip

As one may observe, our interpretation of $\Box _n \varphi$ seems to be very adequate for the subtleties of Brouwer's reasoning. No matter how it came to be that Brouwer stated $R(x) \leftrightarrow \alpha$ (a slip? Lack of notation?) one cannot seriously conclude to KS on the basis of this statement.  So, our conclusion is that $\Box _n \varphi$ in our interpretation is the appropriate instrument for analyzing Brouwer's CS arguments, and in this interpretation TC5 is unattainable, as is KS for which in addition the unattainable decidability is needed.

 \smallskip
 
 We maintain our position that there is no ground for KS  in Brouwer's CS arguments.
 
 \smallskip
 
 In \cite{Niekus2010} we rejected the general Schema of KS but remarked in a footnote the alleged instance of KS above, for which we had, at that moment, no explanation. Possibly the next quote from Mark van Atten's \cite{Atten2014},  is a reaction.

 \begin{smquote}
 John Myhill noticed that KS already had been formulated by Brouwer himself, at least an instance of it. This difference is for some people enough to doubt that Brouwer accepts the general Schema, but I doubt whether one can produce an argument that justifies Brouwer's instance but not the general version \endnote {The whole passage in Dutch: `John Myhill heeft gezien dat Kripke's Schema al door Brouwer zelf geformuleerd is, althans een instantie daarvan. Sommigen zien in dat verschil voldoende grond om te betwijfelen dat Brouwer het algemene schema accepteert, maar er is geloof ik geen argument te bedenken dat wel Brouwer's instantie rechtvaardigt maar niet de algemene versie.Dit omdat de aard en de opbouw van de propositie in kwestie eigenlijk geen rol spelen in het schema zelf. Terecht stellen anderen, waaronder van Dalen, dan ook om voortaan te spreken over ` het Brouwer-Kripke Schema', niet KS maar BKS. }.

\end{smquote}
Our answer is that we have in the above exactly provided the argument van Atten is looking for. Brouwer may have formulated an instance of KS, but he evidently does not want to use it and neither does he mention nor use it anywhere else. Van Atten continues by stating that it would be proper from now on to speak of Brouwer-Kripke Schema instead of Kripke's Schema,  not KS but BKS. And so he does in his `The Creating Subject, the Brouwer-Kripke Schema , and infinite proofs', a defense of the standard theory TCS and KS (\cite{Atten2018}). But the example he uses as his main argument, carries no weight.

\section{Markov's Principle MP}

In the second Vienna Lecture \cite{Brouwer1930A}, Brouwer remarks after each result counting for the full continuum that this result also holds for the full continuum, see section 9. Brouwer does not give any further argument in the Vienna Lecture, but he gives proofs for these results in  \cite{Brouwer1981}, which is the text of Brouwer’s lectures in Cambridge in the period 1946-1951. In this section we shall present such a proof, and we will show that by the present day knowledge concerning Markov’s Principle these  proofs have to be considered as questionable.
 
 In the Vienna lecture Brouwer proves the results for the reduced as well as for the full continuum with a fleeing property $f$. We recall: a fleeing property is  a decidable property for natural numbers; and if $\lambda _f$ is the smallest number possessing $f$, neither $\exists n (\lambda _f= n)$, nor $\neg \exists n (\lambda_f = n)$ has been established. In these proofs Brouwer does not take into consideration whether $\neg \neg \exists n( \lambda _f = n)$ has been established, allthough, as we remarked, it is relevant. In \cite{Brouwer1981} that is different.

 In the Cambridge Lectures (\cite{Brouwer1981}) Brouwer elaborates and refines his analysis of the continuum of the Vienna lecture. For the properties of the reduced continuum he proved in the Vienna Lecture, he now uses a fleeing property for which $\neg \neg \exists n(|lambda _f = n)$ has not been been established.  Each time if it is relevant, Brouwer adds   `$\neg \neg \exists n (\lambda_f = n)$ does not hold' as a condition on $f$. For the properties he proved in Vienna for the full continuum he now uses a choice sequence defined with an $\alpha$  `for which no method of testing is known'.  But he also proves these results for the reduced continuum using a fleeing property for which $\neg \neg \exists n (\lambda_f = n)$ has been established. We shall present such a proof that the reduced continuum is not dense in it self below; the proof for the full continuum we treated in section 9. First we repeat the definition.

\smallskip

A species $S$ is \emph{dense in itself} (in sich dicht) if every point of $S$ is a limiting point.
\smallskip

An element $a$ of a species $S$ is a \emph{limiting point} (Grenzelement) if for a sequence $a_1, a_2, ...$ with $a_1 \succ a_2 \succ .. \succ a$  for every $b  \succ a$ an $a_n$ exists such that $a_n \prec b $, or analogously for a increasing sequence.

\smallskip

Let $a_1, a_2,...$ be a decreasing sequence converging to $0$, let $f$ be a fleeing property, and let $\lambda_f$ its critical number,  such that $\neg \neg \exists n(\lambda_f=n)$ has been established. Define $c_1, c_2,..$ such that $c_i=a_i$  if $i<\lambda_f$, and $c_j= a_{\lambda_f}$ for $j\geq \lambda_f$. The sequence $c_1, c_2, ...$ converges, say to $c$. Then $c\succ 0$, but there is no $n_0$ such that  $a_{n_0} \prec c$ (\cite{Brouwer1981}, p.64).

\bigskip

As usual Brouwer gives no further argument. Our reconstruction, with $\exists n (\lambda _n=n)$ expressed by $\alpha$:
\smallskip

From the definition follows $\exists n \Box _n \alpha \leftrightarrow \exists v (c=a_v)$.

\smallskip
Suppose $c=0$.

\smallskip
Because $c=0 \rightarrow \neg \exists v (c=a_v)$

\smallskip

we would have $\neg \exists n \Box_n \alpha$, and so $\neg \alpha$ (by MD).

\smallskip
Since we have $\neg \neg \exists n (\lambda _f=n)$, i.e. $\neg \neg \alpha$, a contradiction,

\smallskip
so $c \neq 0$.

\smallskip
Since by definition $\neg c<0$ we have $c \succ a$.

\bigskip
Suppose for some $n_0$ that $a_{n_0} \prec c$. 

\smallskip
In that case $\lambda_f \leq n_0$, contradicting that $f$ is fleeing, so we cannot indicate such an $n_0$.

 \bigskip

Brouwer does not give an example of such a fleeing property for which $\neg \neg \exists n (\lambda_f = n)$ has been established. Actually; it is reasonable to assume that he could not have given one, his argument seems to be caught up by later developments. Because if such an example could be formalized in one of the formal systems that exists today, like Heyting's arithmetic HA, it would contradict the validity of \emph{Markov's Rule} which has been established for these systems (\cite{TroelstraDalen1988}, p.136-137). In case of HA it states

\smallskip

MR~~~$\vdash _{HA} \neg \neg \exists x A(x) ~~~  \Rightarrow ~~~ \vdash _{HA} \exists x A(x) $,

\smallskip

Markov's Rule is related to \emph{Markov's Principle}, which states for a decidable property $A$ of natural numbers :

\smallskip

MP~~~ $\neg \neg \exists n A(n) \rightarrow \exists nA(n)$

\smallskip

It has been shown that for HA (as for similar systems) that  MP fails, i.e. a decidable $A$ can be given for which MP is not provable (\cite{Troelstra1993}, pp.~93-94). In fact, in HA it is even so that MP only holds if either $\exists n A(n)$ or $\neg \exists n A(n)$ is provable in HA.

\smallskip

We conclude Brouwer's attempt to show that certain properties for the full continuum can also be proved for the reduced continuum is not satisfactory. We have not been able to establish Brouwer's results using the existence of a fleeing property in another way. So it leaves open the question whether it is possible to prove these properties for the reduced continuum.

\section{Epilogue}

In his Post-War papers Brouwer did not use a special expression for what we have called in this paper a choice sequence. In his Post-War papers he denoted lawlike sequences with `sharp' and all sequences, lawlike or not, by `points'. But in \cite{Brouwer1930A}, his first publication with the method discussed in this  paper, he denoted lawlike sequences with `fertig' and choice sequences with `unfertig', which is German for `complete' and `incomplete`. We used in this paper `choice sequence' to make the connection with the work of Troelstra and van Dalen. But in our \cite{Niekus2010} we used  `incomplete objects' which may be more appropriate. 

In Brouwer's definition a choice sequence is given by a description of a construction which does not determine the values of the sequence completely. In the examples he uses the values of the defined sequence are made dependent of the future mathematical experience of the constructor. The reasoning about a choice sequence is done on the basis of the incomplete definition only, before the construction has started. What holds for a choice sequence is what hold for all of its possible completions.

This conception does not occur in the intuitionistic theories on choice sequences, i.e. the traditional theories and the TCS. For the traditional theories this is easily explained. They are based on a period in which Brouwer did not have the concept himself yet. But if Troelstra states in his final lectures that choice sequences have no mathematical value, only philosophical, he is only speaking of the sequences that exist in his formal systems, and not of the choice sequences Brouwer actually uses.

In so far the TCS is concerned, this theory is based on a misinterpretation. When Kreisel presented the first formalization of the method discussed here, he started with the words:

\begin{smquote}
For general background see the last chapters of Heyting´s book on intuitionism [ \cite{Heyting1956}] and Kleene’s recent monograph (\cite{KleeneVesley1965}). N.B. Knowledge of these two works is a minimal requirement for discussing profitably the present topic; [...]. It is quite unreasonable to want an explanation of intuitionistic notions for `the man in the street´: he simply does not use concepts to which intuitionistic distinctions apply ... (\cite{Kreisel1967}, p.158).
\end{smquote}

`The present topic´ in the quote above is Brouwer´s use of choice sequences as treated in this paper, and they  are perfectly understandable for `the man in the street'. Incomplete objects are all around us in daily life. Most of the measurements yield lawless sequences, i.e. they are given by a finite segment which can be further refined. And sequences in climate prediction, weather reports,  etc. are all choice sequences. It is the acceptance of the incomplete sequences as full-fetched real numbers that makes Brouwer's intuitionism special.

Kreisel's formalization is based on a wrong conception. Although Kreisel´s interpretation has been  changed in the course of time in the intuitionistic literature, from Kreisel's thinking subjects to Troelstra's  IM, from Troelstra's lawlike sequence to Dummett's and van Atten's  choice sequence, there has only been discussion whether C3 should be $\phi \rightarrow \exists n \Box_n \phi$ or $\phi \rightarrow \neg \neg \exists n \Box_n \phi$. The axioms we have rejected, which are $\Box_n \phi \vee \neg \Box_n \phi$ and $\exists n \Box_ \phi \rightarrow \phi$  have never been a subject of discussion. The axiomatization, with the mathematics that has been built on it, has turned out to be a stagnant factor in understanding Brouwer. 

An important conclusion we draw from our conception is that decidability for choice sequences does not hold. For the choice sequence $s$ defined in section 6, and for an arbitrarily $n$, $s_n=0 \vee s_n \neq 0$ does not hold, because we do not know now, whether we shall have made a result in working on the problem involved in the definition. Since the value of $s_n$ is determined at stage n, we do have $\Box_n (s_n=0 \vee s_n \neq 0)$, so $\exists n \Box_n \phi \rightarrow \phi$ is not always valid. This rejection makes it possible to explain Brouwer's subtle distinctions and his proofs, see section 9 and 10, for which the standard theory has no answer. 

The argument for our rejection of decidability is nothing but Aristotle's Sea-battle Argument. In a famous passage of the \emph{Interpretatione} Aristotle  discusses whether the sentence `Tomorrow there will be a sea-battle' has a definite truth value. He concludes that these \emph{future contingents} do not necessarilly have a definite truth value. With our notation we can express Aristotle's hesitation, with $S_m$ denoting `Tomorrow there will be a sea battle', by $S_m \vee \neg S_m$, according to us not true now, and $\Box_m (S_m \vee \neg S_m)$, which is valid.

Finally, our rejection of decidability is fully supported by Brouwer's remarks in the Geneva Lecture, quoted in section 6.  It is the only occasion where Brouwer compares explicitly a particular choice sequence with a lawlike sequence. For the lawlike $r$ Brouwer states that everyone constructing this sequence will reach the same values. For the choice sequence $s$ this does not have to be the case. In our view decidability marks the difference between a choice sequence and a lawlike sequence.

The principles we derived from our conception in section 8 form the core of the Theory of Individual Choice sequences, the TICS. From the material we have presented in this paper we cannot but conclude  that the TICS is far superior to the TCS, both in so far it concerns the underlying conception, as well as for the reconstruction of Brouwer's arguments. But we think our work is of more importance for intuitionism than providing an alternative view on a controversial subject.

Choice sequences are unique for Brouwer's intuitionism. Extended research has been done on choice sequences, resulting in formal systems for which considerable space has been reserved in the handbooks.  In all these systems decidability is presupposed. We claim to have shown on the basis of the choice sequences Brouwer actually uses that decidability is against his ideas.  It seems quite a challenge to incorporate our conclusions concerning undecidability  in the existing systems, but we think it is inescapable for a description of the continuum as Brouwer saw it.

	\theendnotes
\bigskip
    
{Acknowledgement}

I am much indebted to  Dick de Jongh for always being

there for careful reading, comments and support.

	\bibliographystyle{plainnat}

	\bibliography{bibJN}

\end{document}